\documentclass[10pt,noindent,nocenter,sfbold,crosshair,a4paper]{thesis}
\chapapp{lecture}
\usepackage{xy,amssymb,amsmath,diagrams,times,theorem,a4}
\xyoption{curve}
\xyoption{arrow}
\xyoption{matrix}
\usepackage{pstcol,pst-node}
\usepackage{psboxit}

\newcommand{\w}[1]{{\text{\em \usefont{OT1}{cmtt}{m}{n} #1}}}
\newcommand{\wis}[1]{{\text{\em \usefont{OT1}{cmtt}{m}{n} #1}}}
\newcommand{\C}{\mathbb{C}}

\newcommand{\N}{\mathbb{N}}

\newcommand{\Z}{\mathbb{Z}}
\newcommand{\Af}{\mathbb{A}}
\newcommand{\V}{\mathbb{V}}

\newcommand{\Oscr}{\mathcal{O}}
\newcommand{\Ascr}{\mathcal{A}}
\newcommand{\Pscr}{\mathcal{P}}

\newcommand{\vtx}[1]{\boxed{#1}}

\newcommand{\hR}{\hat{R}_{\mathfrak{m}}}

\theoremstyle{plain}
\newtheorem{theorem}{Theorem}[chapter]
\newtheorem{proposition}[theorem]{Proposition}
\newtheorem{lemma}[theorem]{Lemma}
{\theorembodyfont{\rmfamily}
\newtheorem{definition}[theorem]{Definition}

\newtheorem{mexample}[theorem]{Example}
}

 \newcommand{\proofhead}[1]{\par\pagebreak[1]\noindent{\em#1.\ }}
 \newcommand{\pf}{\proofhead{Proof}}
  \newcommand{\qed}{{\unskip\nolinebreak[1]\hspace{1.5em}\mbox{}\nolinebreak
    \hfill$\square$\parfillskip=0pt\finalhyphendemerits=0\par\pagebreak[1]}}

\newenvironment{proof}{\pf}{\qed \par \vskip 6mm}

\begin{document}

\title{three talks on \\ noncommutative geometry@n}

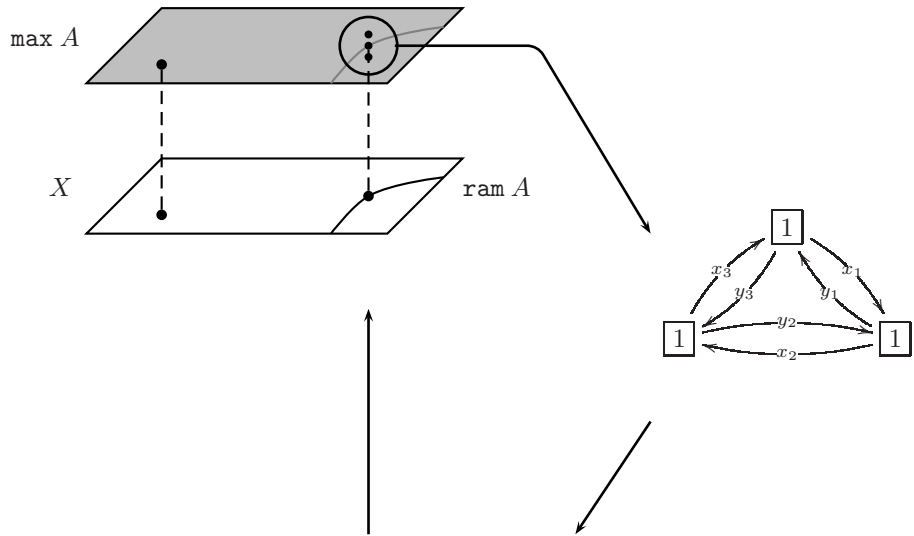
\translator{
\[
\begin{pspicture}(-2,0)(7.5,3)
\psset{xunit=.5cm}
\psset{yunit=.5cm}
\pspolygon[fillstyle=solid,fillcolor=lightgray](0,4)(2,6)(10,6)(8,4)
\pspolygon[fillstyle=solid,fillcolor=white](0,0)(2,2)(10,2)(8,0)
\pscurve[linecolor=gray](6.5,4)(7.5,5)(9.5,5.5)
\pscurve(6.5,0)(7.5,1)(9.5,1.5)
\put(5,.5){$\wis{ram}~A$}
\put(-.5,.5){$X$}
\put(-1,2.5){$\wis{max}~A$}
\cnode*(2,4.5){2pt}{A}
\cnode*(2,0.5){2pt}{B}
\cnode*(7.5,5){1.5pt}{C}
\cnode*(7.5,5.3){1.5pt}{C1}\cnode*(7.5,4.7){1.5pt}{C2}
\pscircle[linewidth=1pt](7.5,5){.4}
\cnode*(7.5,1){2pt}{D}
\ncline[linestyle=dashed]{A}{B}\ncline[linestyle=dashed]{C}{D}
\psline[linewidth=1pt,linearc=.25]{->}(8.2,5)(12,5)(15,0)
\psline[linewidth=1pt,linearc=.25]{->}(15,-5)(13,-8)
\psline[linewidth=1pt,linearc=.25]{->}(7.5,-8)(7.5,-2)
\end{pspicture}~
\xymatrix@=.3cm{
& & \vtx{1} \ar@/^/[lldd]|{y_3} \ar@/^/[rrdd]|{x_1} &&  \\
& & & & \\
\vtx{1} \ar@/^/[rruu]|{x_3} \ar@/^/[rrrr]|{y_2} & & & & \vtx{1} \ar@/^/[lluu]|{y_1} \ar@/^/[llll]|{x_2} }
\]
\vskip 2cm
\[
\hat{A}_{\mathfrak{m}} \simeq \begin{bmatrix} \hR & \hR y_1 + \hR x_3x_2 & \hR x_3+\hR y_1y_2 \\ 
\hR x_1+\hR y_1y_3 & \hR & \hR y_2+\hR x_1x_3 \\
\hR y_3+\hR x_2x_1 & \hR x_2+\hR y_3y_1 & \hR \end{bmatrix} \]
\[
\hR = \C[[x_1y_1,x_2y_2,x_3y_3,x_1x_2x_3,y_1y_2y_3]] \]
}

\author{{\bf lieven le bruyn}}
\institution{university of antwerp}
\dedication{}
\uppertitleback{
}
\middletitleback{}
\lowertitleback{}
\date{2003}
\sloppy

\maketitle

\tableofcontents 

\chapter*{Introduction}

Ever since the dawn of non-commutative algebraic geometry in the mid seventies, see for example the work of 
P. Cohn \cite{Cohn}, J. Golan \cite{Golan}, C. Procesi \cite{Procesibook}, F. Van Oystaeyen and A. Verschoren
\cite{FVO444},\cite{LNM887}, it has been ringtheorists' wet dream that this theory may one day be relevant to
commutative geometry, in particular to the study of singularities and their resolutions.

Over the last decade, non-commutative algebras have been used to construct canonical (partial) resolutions of quotient singularities. That is, take a finite group $G$ acting on $\C^d$ freely away from the origin then the orbit-space $\C^d/G$ is an isolated singularity. Resolutions $Y \rOnto \C^d/G$ have been constructed using the skew group algebra
\[
\C[x_1,\hdots,x_d] \# G \]
which is an order with center $\C[\C^d/G]=\C[x_1,\hdots,x_d]^G$ or deformations of it. In dimension $d=2$ (the case of Kleinian singulariies) this gives us minimal resolutions via the connection with the preprojective algebra, see for example \cite{CrawleyLectNotes}. In dimension $d=3$, the skew group algebra appears via the superpotential and commuting matrices setting (in the physics literature) or via the McKay quiver, see for example \cite{CrawNotes}. If $G$ is Abelian one obtains from this study crepant resolutions but for general $G$ one obtains at best partial resolutions with conifold singularities remaining. In dimension $d > 3$ the situation is unclear at this moment. Usually, skew group algebras and their deformations are studied via homological methods as they are Regular orders, see for example \cite{VdBcrepant}. Here, we will follow a different approach.

My motivation was to find a non-commutative explanation for the omnipresence of conifold singularities in partial resolutions of three-dimensional quotient singularities. Of course you may argue that they have to appear because they are somehow the nicest singularities. But then, what is the corresponding list of 'nice' singularities in dimension four ? or five, six... ?? If my conjectural explanation has any merit the nicest partial resolutions of $\C^4/G$ should contain only singularities which are either polynomials over the conifold or one of the following three types
\[
\frac{\C[[a,b,c,d,e,f]]}{(ae-bd,af-cd,bf-ce)} \qquad \frac{\C[[a,b,c,d,e]]}{(abc-de)} \qquad \frac{\C[[a,b,c,d,e,f,g,h]]}{I} \]
where $I$ is the ideal of all $2 \times 2$ minors of the matrix
\[
\begin{bmatrix}
a & b & c & d \\ e & f & g & h \end{bmatrix}
\]
In dimension $d=5$ the conjecture is that another list of ten new specific singularities will appear, in dimension $d=6$ another $63$ new ones appear and so on.

How do we come to these outlandish conjectures and specific lists? The hope is that any quotient singularity $X = \C^d/G$ has associated to it a 'nice' order $A$ with center $R = \C[X]$ such that there is a stability structure $\theta$ with the scheme of all $\theta$-semistable representations of $A$ being a smooth variety (all these terms will be explained in the main text). If this is the case, the associated moduli space will be a partial resolution
\[
\wis{moduli}^{\theta}_{\alpha}~A \rOnto X = \C^d/G \]
and has a sheaf of Smooth orders $\Ascr$ over it, allowing us to control its singularities in a combinatorial way as depicted in the frontispiece.

If $A$ is a Smooth order over $R = \C[X]$ then its non-commutative variety $\wis{max}~A$ of maximal twosided ideals is birational to $X$ away from the ramification locus. If $P$ is a point of the ramification locus $\wis{ram}~A$ then there is a finite cluster of infinitesimally close non-commutative points lying over it. The local structure of the non-commutative variety $\wis{max}~A$ near this cluster can be summarized by a (marked) quiver setting $(Q,\alpha)$ which in turn allows us to compute the \'etale local structure of $A$ and $R$ in $P$. The central singularities which appear in this way have been classified in \cite{RBLBVdW} up to smooth equivalence giving us the small lists of conjectured singularities.

In these talks I have tried to include background information which may or may not be useful to you. I suggest to browse through the notes by reading the 'jotter-notes' (grey-shaded). If the remark seems obvious to you, carry on. If it puzzles you this may be a good point to enter the main text. More information can be found in the never-ending bookproject \cite{LBnag@n}.

\chapter*{Acknowledgement}

These notes are (hopefully) a streamlined  version of three talks given at the workshop "Sch\'emas de Hilbert, alg\`ebre noncommutative et correspondance de McKay" held at CIRM in Luminy (France),
october 27-31, 2003.

I like to thank the organizers, Jacques Alev, Bernhard Keller and Thierry Levasseur for the invitation (and the possibility to have a nice vacation with part of my family) and the participants for their patience.

\chapter{non-commutative algebra}

The organizers of this conference on "Hilbert schemes, non-commutative algebra and the McKay correspondence" are perfectly aware of my ignorance on Hilbert schemes and  McKay correspondence. I therefore have to assume that I was hired in to tell you something about non-commutative algebra and that is precisely what I intend to do in these three talks.

\section{Why non-commutative algebra?}

Let me begin by trying to motivate why you might get interested in non-commutative algebra if you want to understand quotient singularities and their resolutions. 

So let us take a setting which will be popular this week :  we have a finite group $G$ acting on $d$-dimensional affine space $\C^d$ and this action is  free away from the origin. Then the orbit-space, the so called {\em quotient singularity} $\C^d/G$, is an isolated singularity
\[
\begin{diagram}
\C^d & & \\
\dOnto & & \\
\C^d/G & \lOnto^{res} & Y
\end{diagram}
\]
and we want to construct 'minimal' or 'canonical' resolutions of this singularity. The buzz-word seems to be 'crepant' in these circles. In his Bourbaki talk \cite{ReidTalk} Miles Reid asserts that McKay correspondence follows from a much more general principle

\par \vskip 3mm
\noindent
{\bf Miles Reid's Principle : } Let $M$ be an algebraic manifold, $G$ a group of automorphisms of $M$,
and $Y \rOnto X$ a resolution of singularities of $X = M/G$. Then the answer to any well posed
question about the geometry of $Y$ is the $G$-equivariant geometry of $M$.

\par \vskip 3mm
Applied to the case of quotient singularities, the content of his slogan is that the $G$-equivariant geometry of $\C^d$ already {\em knows} about the crepant resolution $Y \rOnto \C^d/G$. 

Men having principles are an easy target for abuse. So let us change this principle slightly : assume we have an affine variety $M$ on which a reductive group (and for definiteness take $PGL_n$) acts with algebraic quotient variety $M//PGL_n \simeq \C^d/G$
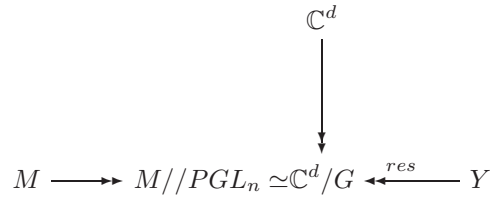
\[
\begin{diagram}
& & & \C^d & & \\
& & & \dOnto & & \\
M & \rOnto & M//PGL_n \simeq & \C^d/G & \lOnto^{res} & Y
\end{diagram}
\]
then, in favorable situations, we can argue that the $PGL_n$-equivariant geometry of $M$ knows about good resolutions $Y$. This brings us to our first entry in our

\begin{boxitpara}{box 0.85 setgray fill}
{\bf jotter  : }

One of the key lessons to be learned from this talk is that $PGL_n$-equivariant geometry of $M$ is {\em roughly} equivalent to the study of a certain non-commutative algebra over $\C^d/G$. In fact, an {\em order} in a central simple algebra of dimension $n^2$ over the function field of the quotient singularity.

Hence, if we know of {\em good} orders over $\C^d/G$, we might get our hands at 'good' resolutions $Y$ by non-commutative methods.
\end{boxitpara}

\section{What non-commutative algebras?}

For the duration of these talks, we will work in the following, quite general, setting : 
\begin{itemize}
\item{$X$ will be a {\em normal} affine variety, possibly having singularities.}
\item{$R$ will be the coordinate ring $\C[X]$ of $X$.}
\item{$K$ will be the function field $\C(X)$ of $X$.}
\end{itemize}
 If you are only interested in quotient singularities, you should replace $X$ by $\C^d/G$, $R$ by the invariant ring $\C[x_1,\hdots,x_d]^G$ and $K$ by the invariant field $\C(x_1,\hdots,x_d)^G$ in all statements below. 

If you are an algebraist, you have my sympathy and our goal will be to construct lots of $R$-orders $A$ in a central simple $K$-algebra $\Sigma$.
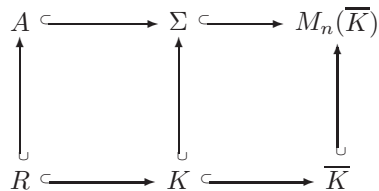
\[
\begin{diagram}
A & \rInto & \Sigma & \rInto & M_n(\overline{K}) \\
\uInto & & \uInto & & \uInto \\
R & \rInto & K & \rInto & \overline{K}
\end{diagram}
\]
If you do not know what a {\em central simple algebra} is, take any non-commutative $K$-algebra $\Sigma$ with center $Z(\Sigma) =K$ such that  over the algebraic closure $\overline{K}$ of $K$ we obtain full $n \times n$ matrices
\[
\Sigma \otimes_K \overline{K} \simeq M_n(\overline{K}) \]
There are plenty such central simple $K$-algebras :

\begin{mexample} \label{cyclic} For any non-zero functions $f,g \in K^*$, the {\em cyclic algebra}
\[
\Sigma = (f,g)_n \qquad \text{defined by } \qquad (f,g)_n = \frac{K \langle x,y \rangle}{(x^n-f,y^n-g,yx-q xy)} \]
with $q$ is a primitive $n$-th root of unity, is a central simple $K$-algebra of dimension $n^2$. Often, $(f,g)_n$ will even be a {\em division algebra}, that is a non-commutative algebra such that every non-zero element has an inverse. 

For example, this is always the case when $E = K[x]$ is a (commutative) field extension of dimension $n$ and if $g$ has order $n$ in the quotient
$K^*/N_{E/K}(E^*)$ where $N_{E/K}$ is the {\em norm map} of $E/K$. See for example \cite[Chp. 15]{Pierce} for more details, but if your German is \'a point I strongly suggest you to read Ina Kersten's book \cite{Kersten} instead.
\end{mexample}

Now, fix such a central simple $K$-algebra $\Sigma$. An {\em $R$-order $A$ in $\Sigma$}  is a subalgebras $A \subset \Sigma$ with center $Z(A) = R$ such that $A$ is finitely generated as an $R$-module and contains a $K$-basis of $\Sigma$, that is
\[
A \otimes_R K \simeq \Sigma \]
The classic reference for orders is Irving Reiner's book \cite{Reiner} but it is hopelessly outdated and focusses  too much on the one-dimensional case. Here is a gap in the market for someone to fill...

\begin{mexample} In the case of quotient singularities $X = \C^d/G$ a natural choice of $R$-order might be the {\em skew group ring} : $\C[x_1,\hdots,x_d] \# G$ which consists of all formal sums $\sum_{g \in G} r_g \# g$ with multiplication defined by
\[
(r \# g)(r' \# g') = r \phi_g(r') \# gg' \]
where $\phi_g$ is the action of $g$ on $\C[x_1,\hdots,x_d]$. The center of the skew group algebra is easily verified to be the ring of $G$-invariants
\[
R = \C[\C^d/G] = \C[x_1,\hdots,x_d]^G \]
Further, one can show that $\C[x_1,\hdots,x_d] \# G$ is an $R$-order in $M_n(K)$ with $n$ the order of $G$. If we ever get to the third lecture, we will give another description of the skew group algebra in terms of the McKay-quiver setting and the variety of commuting matrices.
\end{mexample}

However, there are plenty of other $R$-orders in $M_n(K)$ which may or may not be relevant in the study of the quotient singularity $\C^d/G$.

\begin{mexample} If $f,g \in R - \{ 0 \}$, then the free $R$-submodule of rank $n^2$ of the cyclic $K$-algebra $\Sigma = (f,g)_n$ of example~\ref{cyclic}
\[
A = \sum_{i,j=0}^{n-1} R x^iy^j  \]
is an $R$-order. But there is really no need to go for this 'canonical' example. Someone more twisted may take $I$ and $J$ any two non-zero ideals of $R$, and consider
\[
A_{IJ} = \sum_{i,j=0}^{n-1} I^iJ^j x^iy^j \]
which is an $R$-order too in $\Sigma$ and which is far from being a projective $R$-module unless $I$ and $J$ are invertible $R$-ideals. 

For example, in $M_n(K)$ we can take the 'obvious' $R$-order $M_n(R)$ but one might also take the subring
\[
\begin{bmatrix} R & I \\ J & R \end{bmatrix}
\]
which is an $R$-order if $I$ and $J$ are non-zero ideals of $R$.
\end{mexample}

If you are a geometer (and frankly we are all wannabe geometers these days), our goal is to construct lots of affine $PGL_n$-varieties $M$ such that the algebraic quotient $M // PGL_n$  is  isomorphic to $X$ and, moreover, such that there is a Zariski open subset $U \subset X$
\[
\begin{diagram}
& M & \lInto & \pi^{-1}(U) \\
& \dOnto^{\pi} & & \dOnto_{\text{principal $PGL_n$-fibration}} \\
X \simeq & M//PGL_n & \lInto & U 
\end{diagram}
\]
for which the quotient map is a principal $PGL_n$-fibration, that is, all fibers $\pi^{-1}(u) \simeq PGL_n$ for $u \in U$. 

The connection between such varieties $M$ and orders $A$ in central simple algebras may not be clear at first sight. To give you at least an idea that there is a link,  think of $M$ as the affine variety of $n$-dimensional representations $\w{rep}_n~A$ and of $U$ as the Zariski open subset of all simple $n$-dimensional representations.

Naturally, one can only expect the $R$-order $A$ (or the corresponding $PGL_n$-variety $M$) to be useful in the study of resolutions of $X$ if $A$ is {\em smooth} in some appropriate non-commutative sense. 

Now, there are many  characterizations of {\em commutative} regular domains $R$ : 
\begin{itemize}
\item{$R$ is {\em regular}, that is, has finite global dimension}
\item{$R$ is {\em smooth}, that is, $X$ is a smooth variety}
\end{itemize}
and generalizing either of them to the non-commutative world leads to quite different concepts. 

We will call an $R$-order $A$ is a central simple $K$-algebra $\Sigma$ :
\begin{itemize}
\item{{\bf Regular} if $A$ has finite global dimension together with some extra features such as Auslander regularity or Cohen-Macaulay property, see for example \cite{Levasseur}.}
\item{{\bf Smooth} if the corresponding $PGL_n$-affine variety $M$ is a smooth variety as we will clarify later in this talk.}
\end{itemize}
For applications of Regular orders to desingularizations we refer to the talks by Michel Van den Bergh at this conference or to his paper \cite{VdBcrepant} on this topic. I will concentrate on the properties of Smooth orders instead. Still, it is worth pointing out the strengths and weaknesses of both definitions right now

\begin{boxitpara}{box 0.85 setgray fill}
{\bf jotter  : }

Regular orders are excellent if you want to control homological properties, for example if you want to study the derived categories of their modules. At this moment there is no local characterization of Regular orders if $\wis{dim} X \geq 2$.

 Smooth orders are excellent if you want to have smooth moduli spaces of semi-stable representations. As we will see later, in each dimension there are only a finite number of local types of Smooth orders and these are classified. The downside of this is that Smooth orders are less versatile as Regular orders.
 
 In applications to canonical desingularizations, one often needs the good properties of both so there is a case for investigating SmoothRegular orders better than has been done in the past. 
 \end{boxitpara}
 
In general though, both theories are quite different.

\begin{mexample} The skew group algebra $\C[x_1,\hdots,x_d] \# G$ is always a Regular order but we will see in the next lecture, it is virtually never a Smooth order.
\end{mexample}

\begin{mexample} Let $X$ be the variety of matrix-invariants, that is
\[
X = M_n(\C) \oplus M_n(\C) // PGL_n \]
where $PGL_n$ acts on pairs of $n \times n$ matrices by simultaneous conjugation. The {\em trace ring} of two generic $n \times n$ matrices $A$ is the subalgebra of $M_n(\C[M_n(\C) \oplus M_n(\C)])$ generated over $\C[X]$ by the two {\em generic matrices}
\[
X = \begin{bmatrix} x_{11} & \hdots & x_{1n} \\ 
\vdots & & \vdots \\
x_{n1} & \hdots & x_{nn} \end{bmatrix} \quad \text{and} \quad Y = \begin{bmatrix} y_{11} & \hdots & y_{1n} \\ 
\vdots & & \vdots \\
y_{n1} & \hdots & y_{nn} \end{bmatrix} \]
Then, $A$ is an $R$-order in a division algebra of dimension $n^2$ over $K$, called the {\em generic division algebra}. Moreover, $A$ is a Smooth order but is Regular only when $n=2$, see
\cite{LBVdB}.
\end{mexample}

\section{Constructing orders by descent}

\begin{boxitpara}{box 0.85 setgray fill}
{\bf jotter  : }

French mathematicians have developed in the sixties an elegant theory, called {\em descent theory}, which allows one to construct elaborate examples out of trivial ones by bringing in topology. This theory allows to classify objects which are only {\em locally} (but not necessarily globally) trivial.
\end{boxitpara}

For applications to orders there are two topologies to consider : the well-known Zariski topology and the perhaps lesser-known \'etale topology.
Let us try to give a formal definition of Zariski and \'etale {\em covers} aimed at ringtheorists.

A {\em Zariski cover} of $X$ is a finite product of localizations at elements of $R$
\[
S_z = \prod_{i=1}^k R_{f_i} \qquad \text{such that } \qquad (f_1,\hdots,f_k) = R \]
and is therefore a faithfully flat extension of $R$. Geometrically, the ringmorphism $R \rTo S_z$ defines a cover of $X = \w{spec}~R$ by $k$ disjoint sheets $\w{spec}~S_z = \sqcup_i \w{spec}~R_{f_i}$, each corresponding to a Zariski open subset of $X$, the complement of $\V(f_i)$ and the condition is that these closed subsets $\V(f_i)$ do not have a point in common.
That is, we have the picture of figure~1.1 :
\begin{figure} \label{Zariski}
\[
\begin{pspicture}(-1,0)(7.5,6)
\psset{xunit=.5cm}
\psset{yunit=.5cm}
\pspolygon[fillstyle=solid,fillcolor=lightgray](0,10)(2,12)(10,12)(8,10)
\pspolygon[fillstyle=solid,fillcolor=gray](0,8)(2,10)(10,10)(8,8)
\pspolygon[fillstyle=solid,fillcolor=lightgray](0,4)(2,6)(10,6)(8,4)
\pspolygon[fillstyle=solid,fillcolor=white](0,0)(2,2)(10,2)(8,0)
\pscurve[linecolor=gray](0.3,10)(2,11)(4,11.3)(6,12)
\pscurve(.3,0)(2,1)(4,1.3)(6,2)
\pscurve[linecolor=lightgray](1,9)(4,9)(6.5,10)
\pscurve(1,1)(4,1)(6.5,2)
\pscurve[linecolor=gray](6.5,4)(7.5,5)(9.5,5.5)
\pscurve(6.5,0)(7.5,1)(9.5,1.5)
\put(5,.5){$\wis{spec}~R$}
\put(5,2.5){$\wis{spec}~R_{f_k}$}
\put(5,4.5){$\wis{spec}~R_{f_2}$}
\put(5,5.5){$\wis{spec}~R_{f_1}$}
\put(3,3.5){$\vdots$}
\end{pspicture}
\]
\caption{A Zariski cover of $X=\wis{spec}~R$}
\end{figure}

Zariski covers form a {\em Grothendieck topology}, that is, two Zariski covers $S^1_z = \prod_{i=1}^k R_{f_i}$ and $S^2_z = \prod_{j=1}^l R_{g_j}$ have a common refinement
\[
S_z = S^1_z \otimes_R S^2_z = \prod_{i=1}^k \prod_{j=1}^l R_{f_ig_j} \]
For a given Zariski cover $S_z = \prod_{i=1}^k R_{f_i}$ a corresponding {\em \'etale cover} is a product
\[
S_e = \prod_{i=1}^k \frac{R_{f_i}[x(i)_1,\hdots,x(i)_{k_i}]}{(g(i)_1,\hdots,g(i)_{k_i})} \quad \text{with} \quad \begin{bmatrix}
\frac{\partial g(i)_1}{\partial x(i)_1} & \hdots & \frac{\partial g(i)_1}{\partial x(i)_{k_i}} \\
\vdots & & \vdots \\
\frac{\partial g(i)_{k_i}}{\partial x(i)_1} & \hdots & \frac{\partial g(i)_{k_i}}{\partial x(i)_{k_i}} 
\end{bmatrix} \]
a unit in the $i$-th component of $S_e$. In fact, for applications to orders it is usually enough to consider {\em special etale extensions}
\[
S_e = \prod_{i=1}^k \frac{R_{f_i}[x]}{(x^{k_i} - a_i)} \quad \text{where} \quad a_i~\text{is a unit in $R_{f_i}$} \]
Geometrically, an \'etale cover determines for every Zariski sheet $\w{spec}~R_{f_i}$ a {\em locally isomorphic} (for the analytic topology) multi-covering and the number of sheets may vary with $i$ (depending on the degrees of the polynomials $g(i)_j \in R_{f_i}[x(i)_1,\hdots,x(i)_{k_i}]$. That is, the mental picture corresponding to an \'etale cover is given in figure~1.2  below.
\begin{figure} \label{Etale}
\[
\begin{pspicture}(-1,0)(7.5,9)
\psset{xunit=.5cm}
\psset{yunit=.5cm}
\pspolygon[fillstyle=solid,fillcolor=lightgray](0,12)(2,14)(10,14)(8,12)
\pspolygon[fillstyle=solid,fillcolor=lightgray](0,14)(2,16)(10,16)(8,14)
\pspolygon[fillstyle=solid,fillcolor=lightgray](0,16)(2,18)(10,18)(8,16)
\pspolygon[fillstyle=solid,fillcolor=gray](0,8)(2,10)(10,10)(8,8)
\pspolygon[fillstyle=solid,fillcolor=gray](0,10)(2,12)(10,12)(8,10)
\pspolygon[fillstyle=solid,fillcolor=lightgray](0,4)(2,6)(10,6)(8,4)
\pspolygon[fillstyle=solid,fillcolor=white](0,0)(2,2)(10,2)(8,0)
\pscurve[linecolor=gray](0.3,12)(2,13)(4,13.3)(6,14)
\pscurve[linecolor=gray](0.3,14)(2,15)(4,15.3)(6,16)
\pscurve[linecolor=gray](0.3,16)(2,17)(4,17.3)(6,18)
\pscurve(.3,0)(2,1)(4,1.3)(6,2)
\pscurve[linecolor=lightgray](1,9)(4,9)(6.5,10)
\pscurve[linecolor=lightgray](1,11)(4,11)(6.5,12)
\pscurve(1,1)(4,1)(6.5,2)
\pscurve[linecolor=gray](6.5,4)(7.5,5)(9.5,5.5)
\pscurve(6.5,0)(7.5,1)(9.5,1.5)
\put(5,.5){$\wis{spec}~R$}
\put(5,2.5){$\wis{spec}~\frac{R_{f_k}[x(k)_1,\hdots,x(k)_{k_k}]}{(g(k)_1,\hdots,g(k)_{k_k})}$}
\put(5,5.2){$\wis{spec}~\frac{R_{f_2}[x(2)_1,\hdots,x(2)_{k_2}]}{(g(2)_1,\hdots,g(2)_{k_2})}$}
\put(5,7.5){$\wis{spec}~\frac{R_{f_1}[x(1)_1,\hdots,x(1)_{k_1}]}{(g(1)_1,\hdots,g(1)_{k_1})}$}
\put(3,3.5){$\vdots$}
\end{pspicture}
\]
\caption{An \'etale cover of $X = \wis{spec}~R$}
\end{figure}

Again, \'etale covers form a Zariski topology as the common refinement $S^1_e \otimes_R S^2_e$ of two \'etale covers is again \'etale because its components are of the form
\[
\frac{R_{f_ig_j}[x(i)_1,\hdots,x(i)_{k_i},y(j)_1,\hdots,y(j)_{l_j}]}{(g(i)_1,\hdots,g(i)_{k_i},h(j)_1,\hdots,h(j)_{l_j})} \]
and the Jacobian-matrix condition for each of these components is again satisfied. Because of the local isomorphism property many ringtheoretical local properties (such as smoothness, normality etc.) are preserved under \'etale covers.

Now, fix an $R$-order $B$ in some central simple $K$-algebra $\Sigma$, then a {\em Zariski twisted form} $A$ of $B$ is an $R$-algebra such that 
\[
A \otimes_R S_z \simeq B \otimes_R S_z \]
for some Zariski cover $S_z$ of $R$. 

If $P \in X$ is a point with corresponding maximal ideal $\mathfrak{m}$, then $P \in \wis{spec}~R_{f_i}$ for some of the components of $S_z$ and as $A_{f_i} \simeq B_{f_i}$ we have for the local rings at $P$
\[
A_{\mathfrak{m}} \simeq B_{\mathfrak{m}} \]
that is, the Zariski local information of any Zariski-twisted form of $B$ is that of $B$ itself.

Likewise, an {\em \'etale twisted form} $A$ of $B$ is an $R$-algebra such that
\[
A \otimes_R S_e \simeq B \otimes_R S_e \]
for some \'etale cover $S_e$ of $R$. 

This time the Zariski local information of $A$ and $B$ may be different at a point $P \in X$ but we do have that the $\mathfrak{m}$-adic completions of $A$ and $B$
\[
\hat{A}_{\mathfrak{m}} \simeq \hat{B}_{\mathfrak{m}} \]
are isomorphic as $\hat{R}_{\mathfrak{m}}$-algebras. 

\begin{boxitpara}{box 0.85 setgray fill}
{\bf jotter  : }

The Zariski local structure of $A$ determines the localization $A_{\mathfrak{m}}$, the \'etale local structure determines the completion $\hat{A}_{\mathfrak{m}}$.
\end{boxitpara}

Descent theory allows to classify Zariski- or \'etale twisted forms of an $R$-order $B$ by means of the corresponding cohomology groups of the automorphism schemes. For more details on this please read the book \cite{KnusOjanguren} by M. Knus and M. Ojanguren if you are a ringtheorist and that of S. Milne \cite{Milne} if you are more of a geometer. 

If one applies descent to the most trivial of all $R$-orders, the full matrix algebra $M_n(R)$, one arrives at

\section{Azumaya algebras}

A Zariski twisted form of $M_n(R)$ is an $R$-algebra $A$ such that
\[
A \otimes_R S_z \simeq M_n(S_z) = \prod_{i=1}^k M_n(R_{f_i}) \]
Conversely, you can construct such twisted forms by {\em gluing together} the matrix rings $M_n(R_{f_i})$. The easiest way to do this is to glue $M_n(R_{f_i})$ with $M_n(R_{f_j})$ over $R_{f_if_j}$ via the natural embeddings
\[
R_{f_i} \rInto R_{f_if_j} \lInto R_{f_j} \]
Not surprisingly, we obtain in this way $M_n(R)$ back. 

But there are more clever ways to perform the gluing by bringing in the non-commutativity of matrix-rings. We can glue
\[
M_n(R_{f_i})  \rInto  M_n(R_{f_if_j})  \rTo^{g_{ij}.g_{ij}^{-1}}_{\simeq}  M_n(R_{f_if_j})  \lInto M_n(R_{f_j})
\]
over their intersection via conjugation with an invertible matrix $g_{ij} \in GL_n(R_{f_if_j})$. If the elements $g_{ij}$ for $1 \leq i,j \leq k$ satisfy the {\em cocycle condition} (meaning that the different possible gluings are compatible over their common localization $R_{f_if_jf_l}$), we obtain a sheaf of non-commutative algebras $\mathcal{A}$ over $X = \wis{spec}~R$ such that its global sections are not necessarily $M_n(R)$. 

\begin{proposition} Any Zariski twisted form of $M_n(R)$ is isomorphic to
\[
End_R(P) \]
where $P$ is a projective $R$-module of rank $n$. Two such twisted forms are isomorphic as $R$-algebras
\[
End_R(P) \simeq End_R(Q) \quad \text{iff} \quad P \simeq Q \otimes I \]
for some invertible $R$-ideal $I$.
\end{proposition}

\begin{proof}[sketch] We have an exact sequence of groupschemes
\[
1 \rTo \mathbb{G}_m \rTo \wis{GL}_n \rTo \wis{PGL}_n \rTo 1 \]
(here, $\mathbb{G}_m$ is the sheaf of units) and taking Zariski cohomology groups over $X$ we have a sequence
\[
1 \rTo H^1_{Zar}(X,\mathbb{G}_m) \rTo H^1_{Zar}(X,\wis{GL}_n) \rTo H^1_{Zar}(X,\wis{PGL}_n) \]
where the first term is isomorphic to the Picard group $Pic(R)$ and the second term classifies projective $R$-modules of rank $n$ upto isomorphism. The final term classifies the Zariski twisted forms of $M_n(R)$ as the automorphism group of $M_n(R)$ is $PGL_n$.
\end{proof}

\begin{mexample} Let $I$ and $J$ be two invertible ideals of $R$, then
\[
End_R(I \oplus J) \simeq \begin{bmatrix} R & I^{-1}J \\ I J^{-1} & R \end{bmatrix} \subset M_2(K) \]
and if $IJ^{-1} = (r)$ then $I \oplus J \simeq (Rr \oplus R) \otimes J$ and indeed we have an isomorphism
\[
\begin{bmatrix} 1 & 0 \\ 0 & r^{-1} \end{bmatrix} \begin{bmatrix} R & I^{-1} J \\ I J^{-1} & R \end{bmatrix} \begin{bmatrix} 1 & 0 \\ 0 & r  \end{bmatrix} = \begin{bmatrix} R & R \\ R & R \end{bmatrix}
\]
\end{mexample}

Things get a lot more interesting in the \'etale topology.

\begin{definition} An {\em $n$-Azumaya algebra} over $R$ is an \'etale twisted form $A$ of $M_n(R)$. If $A$ is also a Zariski twisted form we call $A$ a {\em trivial} Azumaya algebra.
\end{definition}

From the definition and faithfully flat descent, the following facts follow :

\begin{lemma} \label{propAzu} If $A$ is an $n$-Azumaya algebra over $R$, then :
\begin{enumerate}
\item{The center $Z(A) = R$ and $A$ is a projective $R$-module of rank $n^2$.}
\item{All simple $A$-representations have dimension $n$ and for every maximal ideal $\mathfrak{m}$ of $R$ we have
\[
A / \mathfrak{m} A \simeq M_n(\C) \]}
\end{enumerate}
\end{lemma}

\begin{proof} For (2) take $M \cap R = \mathfrak{m}$ where $M$ is the kernel of a simple representation $A \rOnto M_k(\C)$, then as $\hat{A}_{\mathfrak{m}} \simeq M_n(\hat{R}_{\mathfrak{m}})$ it follows that
\[
A / \mathfrak{m} A \simeq M_n(\C) \]
and hence that $k=n$ and $M = A \mathfrak{m}$.
\end{proof}

It is clear from the definition that when $A$ is an $n$-Azumaya algebra and $A'$ is an $m$-Azumaya algebra over $R$, $A \otimes_R A'$ is an $mn$-Azumaya and also that
\[
A \otimes_R A^{op} \simeq End_R(A) \]
where $A^{op}$ is the {\em opposite} algebra (that is, equipped with the reverse multiplication rule). 

These facts allow us to define the {\em Brauer group} $\wis{Br} R$ to be the set of equivalence classes $[A]$ of Azumaya algebras over $R$ where
\[
[A] = [A'] \qquad \text{iff} \qquad A \otimes_R A' \simeq End_R(P) \]
for some projective $R$-module $P$ and where multiplication is induced from the rule
\[
[A].[A'] = [A \otimes_R A'] \]
One can extend the definition of the Brauer group from affine varieties to arbitrary schemes and A. Grothendieck has shown that the Brauer group of a projective smooth variety is a birational invariant, see \cite{Grothendieck}. Moreover, he conjectured a cohomological description of the Brauer group $\wis{Br} R$ which was subsequently proved by O. Gabber in \cite{Gabber}.

\begin{theorem} The Brauer group is an \'etale cohomology group
\[
\wis{Br} R \simeq H^2_{et}(X, \mathbb{G}_m)_{torsion} \]
where $\mathbb{G}_m$ is the unit sheaf and where the subscript denotes that we take only torsion elements. If $R$ is regular, then $H^2_{et}(X, \mathbb{G}_m)$ is torsion so we can forget the subscript.
\end{theorem}

This result should be viewed as the ringtheory analogon of the {\em crossed product theorem} for central simple algebras over fields, see for example \cite{Pierce}.

Observe that in Gabber's result there is no sign of singularities in the description of the Brauer group. In fact, with respect to the desingularization problem, Azumaya algebras are only as good as their centers.

\begin{proposition} If $A$ is an $n$-Azumaya algebra over $R$, then 
\begin{enumerate}
\item{$A$ is Regular iff $R$ is commutative regular.}
\item{$A$ is Smooth iff $R$ is commutative regular.}
\end{enumerate}
\end{proposition}

\begin{proof}
$(1)$ follows from faithfully flat descent and $(2)$ from lemma~\ref{propAzu} which asserts that the $PGL_n$-affine variety corresponding to $A$ is a principal $PGL_n$-fibration in the \'etale topology, which shows that both $n$-Azumaya algebras and principal $PGL_n$-fibrations are classified by the \'etale cohomology group $H^1_{et}(X,\w{PGL}_n)$.
\end{proof}

\begin{boxitpara}{box 0.85 setgray fill}
{\bf jotter  : }

In the correspondence between $R$-orders and $PGL_n$-varieties, Azumaya algebras correspond to {\em principal} $PGL_n$-fibrations over $X$. With respect to desingularizations, Azumaya algebras are therefore only as good as their centers.
\end{boxitpara}

\section{Reflexive Azumaya algebras}

So let us bring in {\em ramification} in order to construct orders which may be more useful in our desingularization project.

\begin{mexample} \label{exram} Consider the $R$-order in $M_2(K)$
\[
A = \begin{bmatrix} R & R \\ I & R \end{bmatrix} \]
where $I$ is some ideal of $R$. If $P \in X$ is a point with corresponding maximal ideal $\mathfrak{m}$ we have that :

For $I$ not contained in $\mathfrak{m}$ we have $A_{\mathfrak{m}} \simeq M_2(R_{\mathfrak{m}})$ whence $A$ is an Azumaya algebra in $P$.

For $I \subset \mathfrak{m}$ we have
\[
A_{\mathfrak{m}} \simeq \begin{bmatrix} R_{\mathfrak{m}} & R_{\mathfrak{m}} \\ I_{\mathfrak{m}} & R_{\mathfrak{m}} \end{bmatrix} \not= M_2(R_{\mathfrak{m}}) \]
whence $A$ is not Azumaya in $P$.
\end{mexample}

\begin{definition} The {\em ramification locus} of an $R$-order $A$ is the Zariski closed subset of $X$ consisting of those points $P$ such that for the corresponding maximal ideal $\mathfrak{m}$
\[
A/\mathfrak{m}A \not\simeq M_n(\C) \]
That is, $\wis{ram}~A$ is the locus of $X$ where $A$ is not an Azumaya algebra. Its complement $\wis{azu}~A$ is called the {\em Azumaya locus} of $A$ which is always a Zariski open subset of $X$.
\end{definition}

\begin{definition} An $R$-order $A$ is said to be a {\em reflexive $n$-Azumaya algebra} iff
\begin{enumerate}
\item{$\wis{ram}~A$ has codimension at least two in $X$, and}
\item{$A$ is a reflexive $R$-module}
\end{enumerate}
that is, $A \simeq Hom_R(Hom_R(A,R),R) = A^{**}$.
\end{definition}

The origin of the terminology is that when $A$ is a reflexive $n$-Azumaya algebra we have that $A_{\mathfrak{p}}$ is $n$-Azumaya for every height one prime ideal $\mathfrak{p}$ of $R$ and that $A = \cap_{\mathfrak{p}} A_{\mathfrak{p}}$ where the intersection is taken over all height one primes.

 For example, in example~\ref{exram} if $I$ is a divisorial ideal of $R$, then $A$ is not reflexive Azumaya as $A_{\mathfrak{p}}$ is not Azumaya for $\mathfrak{p}$ a height one prime containing $I$ and if $I$ has at least height two, then $A$ is often not a reflexive Azumaya algebra because $A$ is not reflexive as an $R$-module. For example take 
 \[
A =  \begin{bmatrix}
 \C[x,y] & \C[x,y] \\
 (x,y) & \C[x,y] 
 \end{bmatrix} 
 \]
 then the reflexive closure of  $A$ is $A^{**} = M_2(\C[x,y])$.

Sometimes though, we get reflexivity of $A$ for free, for example when $A$ is a Cohen-Macaulay $R$-module. An other important fact to remember is that for $A$ a reflexive Azumaya, $A$ is Azumaya if and only if $A$ is projective as an $R$-module. If you want to know more about reflexive Azumaya algebras you may want to read \cite{Orzech} or my Ph.D. thesis \cite{LBphd}.

\begin{mexample} Let $A = \C[x_1,\hdots,x_d] \# G$ then $A$ is a reflexive Azumaya algebra whenever $G$ acts freely away from the origin and $d \geq 2$. Moreover, $A$ is never an Azumaya algebra as its ramification locus is the isolated singularity.
\end{mexample}

In analogy with the Brauer group one can define the {\em reflexive Brauer group} $\beta(R)$ whose elements are the equivalence classes $[A]$ for $A$ a reflexive Azumaya algebra over $R$ with equivalence relation
\[
[A] = [A'] \qquad \text{iff} \qquad A \otimes_R A' \simeq End_R(M) \]
where $M$ is a reflexive $R$-module and with multiplication induced by the rule
\[
[A].[A'] = [ (A \otimes_R A')^{**} ] \]
In \cite{LBreflexiveAzu} it was shown that the reflexive Brauer group does have a cohomological description

\begin{proposition} The reflexive Brauer group is an \'etale cohomology group
\[
\beta(R) \simeq H^2_{et}(X_{sm},\mathbb{G}_m) \]
where $X_{sm}$ is the smooth locus of $X$.
\end{proposition}

This time we see that the singularities of $X$ do appear in the description so perhaps reflexive Azumaya algebras are a class of orders more suitable for our project. This is even more evident if we impose non-commutative smoothness conditions on $A$.

\begin{proposition} Let $A$ be a reflexive Azumaya algebra over $R$, then :
\begin{enumerate}
\item{if $A$ is Regular, then $\wis{ram}~A = X_{sing}$, and}
\item{if $A$ is Smooth, then $X_{sing}$ is contained in $\wis{ram}~A$.}
\end{enumerate}
\end{proposition}

\begin{proof}
$(1)$ was proved in \cite{LBcentralsing} the essential point being that if $A$ is Regular then $A$ is a Cohen-Macaulay $R$-module whence it must be projective over a smooth point of $X$ but then it is not just an reflexive Azumaya but actually an Azumaya algebra in that point. The second statement can be further refined as we will see in the next lecture.
\end{proof}

Many classes of well-studied algebras are reflexive Azumaya algebras, 
\begin{itemize}
\item{Trace rings $\mathbb{T}_{m,n}$ of $m$ generic $n \times n$ matrices (unless $(m,n)=(2,2)$), see \cite{LBArtinSchofield}.}
\item{Quantum enveloping algebras $U_q(\mathfrak{g})$ of semi-simple Lie algebras at roots of unity, see for example \cite{BG1}.}
\item{Quantum function algebras $O_q(G)$ for semi-simple Lie groups at roots of unity, see for example \cite{BG2}.}
\item{Symplectic reflection algebras $A_{t,c}$, see \cite{BG3}.}
\end{itemize}

\begin{boxitpara}{box 0.85 setgray fill}
{\bf jotter  : }

Many interesting classes of Regular orders are reflexive Azumaya algebras. As a consequence their ramification locus coincides with the singularity locus of the center.
\end{boxitpara}

\section{Cayley-Hamilton algebras}

It is about time to clarify the connection with $PGL_n$-equivariant geometry. We will introduce a class of non-commutative algebras, the so called {\em Cayley-Hamilton algebras} which are the level $n$ generalization of the category of commutative algebras and which contain all $R$-orders.

A {\em trace map} $tr$ is a $\C$-linear function $A \rTo A$ satisfying for all $a,b \in A$
\[
tr(tr(a)b) = tr(a)tr(b) \qquad tr(ab) = tr(ba) \qquad \text{and} \qquad tr(a)b = b tr(a) \]
so in particular, the image $tr(A)$ is contained in the center of $A$.

If $M \in M_n(R)$ where $R$ is a commutative $\C$-algebra, then its characteristic polynomial
\[
\chi_M = det(t1_n-M) = t^n + a_1 t^{n-1} + a_2 t^{n-2} + \hdots + a_n \]
has coefficients $a_i$ which are polynomials with rational coefficients in traces of powers of $M$
\[
a_i = f_i(tr(M),tr(M^2),\hdots,tr(M^{n-1}) \]
Hence, if we have an algebra $A$ with a trace map $tr$ we can define a {\em formal characteristic polynomial} of degree $n$ for every $a \in A$ by taking
\[
\chi_a = t^n + f_1(tr(a),\hdots,tr(a^{n-1}) t^{n-1} + \hdots + f_n(tr(a),\hdots,tr(a^{n-1}) \]
which allows us to define the category $\wis{alg@n}$ of Cayley-Hamilton algebras of degree $n$.

\begin{definition} An object $A$ in $\wis{alg@n}$ is a Cayley-Hamilton algebra of degree $n$, that is, a $\C$-algebra with trace map $tr~:~A \rTo A$ satisfying
\[
tr(1) = n \qquad \text{and} \qquad \forall a \in A~:~\chi_a(a) = 0 \]
Morphisms $A \rTo B$ in $\wis{alg@n}$ are trace preserving $\C$-algebra morphisms, that is,
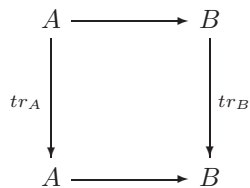
\[
\begin{diagram}
A & \rTo & B \\
\dTo^{tr_A} & & \dTo_{tr_B} \\
A & \rTo & B \end{diagram}
\]
is a commutative diagram.
\end{definition}

\begin{mexample} Azumaya algebras, reflexive Azumaya algebras and more generally every $R$-order $A$ in a central simple $K$-algebra of dimension $n^2$ is a Cayley-Hamilton algebra of degree $n$. For, consider the inclusions
\[
\begin{diagram}
A & \rInto & \Sigma & \rInto & M_n(\overline{K}) \\
\dDotsto^{tr} & & \dDotsto^{tr} & & \dTo^{tr} \\
R & \rInto & K & \rInto & \overline{K}
\end{diagram}
\]
Here, $tr~:~M_n(\overline{K}) \rTo \overline{K}$ is the usual trace map. By Galois descent this induces a trace map, the so called {\em reduced trace}, $tr~:~\Sigma \rTo K$. Finally, because $R$ is integrally closed in $K$ and $A$ is a finitely generated $R$-module it follows that $tr(a) \in R$ for every element $a \in A$.
\end{mexample}

 If $A$ is a finitely generated object in $\wis{alg@n}$, we can define an affine $PGL_n$-scheme,
 $\wis{trep}_n~A$, 
 classifying all trace preserving $n$-dimensional representations $A \rTo^{\phi} M_n(\C)$ of $A$. The action of $PGL_n$ on $\wis{trep}_n~A$ is induced by conjugation in the target space, that is
 $g.\phi$ is the trace preserving algebra map
 \[
 A \rTo^{\phi} M_n(\C) \rTo^{g_g^{-1}} M_n(\C) \]
 Orbits under this action correspond precisely to isomorphism classes of representations. The scheme $\wis{trep}_n~A$ is a closed subscheme of $\wis{rep}_n~A$ the more familiar $PGL_n$-affine scheme of all $n$-dimensional representations of $A$. In general, both schemes may be different.
 
 \begin{mexample} Let $A$ be the quantum plane at $-1$, that is
 \[
 A = \frac{\C \langle x,y \rangle}{(xy+yx)} \]
 then $A$ is an order with center $R=\C[x^2,y^2]$ in the quaternion algebra $(x,y)_2 = K1 \oplus Ku \oplus Kv \oplus Kuv$ over $K = \C(x,y)$ where $u^2=x.v^2=y$ and $uv=-vu$. Observe that $tr(x) = tr(y) = 0$ as the embedding $A \rInto (x,y)_2 \rInto M_2(\C[u,y])$ is given by
 \[
 x \mapsto \begin{bmatrix} u & 0 \\ 0 & -u \end{bmatrix} \qquad \text{and} \qquad
 y \mapsto \begin{bmatrix} 0 & 1 \\ y & 0 \end{bmatrix} \]
 Therefore, a trace preserving algebra map $A \rTo M_2(\C)$ is fully determined by the images of $x$ and $y$ which are trace zero $2 \times 2$ matrices
 \[
 \phi(x) = \begin{bmatrix} a & b \\ c & -a \end{bmatrix} \quad \text{and} \quad \phi(y) = \begin{bmatrix} d & e \\ f & -d \end{bmatrix} \quad \text{satisfying} \quad bf+ce = 0 \]
 That is, $\wis{trep}_2~A$ is the hypersurface $\V(bf+ce) \subset \Af^6$ which has a unique isolated singularity at the origin. However, $\wis{rep}_2~A$ contains more points, for example
 \[
 \phi(x) = \begin{bmatrix} a & 0 \\ 0 & b \end{bmatrix} \quad \text{and} \quad
 \phi(y) = \begin{bmatrix} 0 & 0 \\ 0 & 0 \end{bmatrix} \]
 is a point in $\wis{rep}_2~A - \wis{trep}_2~A$ whenever $b \not= -a$.
 \end{mexample}
 
 A functorial description of $\wis{trep}_n~A$ is given by the following universal property proved by C. Procesi \cite{ProcesiCH}
 
 \begin{theorem} Let $A$ be a $\C$-algebra with trace map $tr_A$, then there is a trace preserving algebra morphism
 \[
 j_A~:~A \rTo M_n(\C[\wis{trep}_n~A]) \]
 satisfying the following universal property. If $C$ is a commutative $\C$-algebra and there is a trace preserving algebra map $A \rTo^{\psi} M_n(C)$ (with the usual trace on $M_n(C)$), then there is a unique algebra morphism $\C[\wis{trep}_n~A] \rTo^{\phi} C$ such that the diagram
 \[
 \begin{diagram}
 A & \rTo^{\psi} & M_n(C) \\
 \dTo^{j_A} & \ruTo^{M_n(\phi)} & \\
 M_n(\C[\wis{trep}_n~A]) & &
 \end{diagram}
 \]
 is commutative. Moreover, $A$ is an object in $\wis{alg@n}$ if and only if $j_A$ is a monomorphism.
 \end{theorem}
 
 The $PGL_n$-action on $\wis{trep}_n~A$ induces an action of $PGL_n$ by automorphisms on $\C[\wis{trep}_n~A]$. On the other hand, $PGL_n$ acts by conjugation on $M_n(\C)$ so we have a combined action on $M_n(\C[\wis{trep}_n~A]) = M_n(\C) \otimes \C[\wis{trep}_n~A]$ and it follows from the universal property that the image of $j_A$ is contained in the ring of $PGL_n$-invariants
 \[
 A \rTo^{j_A} M_n(\C[\wis{trep}_n~A])^{PGL_n} \]
 which is an inclusion if $A$ is a Cayley-Hamilton algebra. In fact, C. Procesi proved in \cite{ProcesiCH} the following important result which allows to reconstruct orders and their centers from $PGL_n$-equivariant geometry.
 
 \begin{theorem} The functor
 \[
 \wis{trep}_n~:~\wis{alg@n} \rTo \wis{PGL(n)-affine} \]
 has a {\em left} inverse
 \[
 \wis{A}_{-}~:~\wis{PGL(n)-affine} \rTo \wis{alg@n} \]
 defined by $\wis{A}_Y = M_n(\C[Y])^{PGL_n}$. In particular, we have for any $A$ in $\wis{alg@n}$
 \[
 A = M_n(\C[\wis{trep}_n~A])^{PGL_n} \qquad \text{and} \qquad tr(A) = \C[\wis{trep}_n~A]^{PGL_n}
 \]
 That is the central subalgebra $tr(A)$ is the coordinate ring of the algebraic quotient variety
 \[
 \wis{trep}_n~A // PGL_n = \wis{tiss}_n~A \]
 classifying isomorphism classes of trace preserving semi-simple $n$-dimensional representations of $A$.
 \end{theorem}
 
However, these functors do {\em not} give an equivalence between $\wis{alg@n}$ and $PGL_n$-equivariant affine geometry. There are plenty more $PGL_n$-varieties than Cayley-Hamilton algebras.

\begin{mexample} Conjugacy classes of nilpotent matrices in $M_n(\C)$ correspond bijective to partitions $\lambda = (\lambda_1 \geq \lambda_2 \geq \hdots )$ of $n$ (the $\lambda_i$ determine the sizes of the Jordan blocks). It follows from the Gerstenhaber-Hesselink theorem that the closures of such orbits
\[
\overline{\Oscr_{\lambda}} = \cup_{\mu \leq \lambda} \Oscr_{\mu} \]
where $\leq$ is the dominance order relation. Each $\overline{\Oscr_{\lambda}}$ is an affine $PGL_n$-variety and the corresponding algebra is
\[
\wis{A}_{\overline{\Oscr_{\lambda}}} = \C[x]/(x^{\lambda_1}) \]
whence many orbit closures (all of which are affine $PGL_n$-varieties) correspond to the same algebra.
\end{mexample}

\begin{boxitpara}{box 0.85 setgray fill}
{\bf jotter  : }

The category $\wis{alg@n}$ of Cayley-Hamilton algebras is to $\wis{noncommutative geometry@n}$ what $\wis{commalg}$, the category of all commutative algebras is to commutative algebraic geometry.

In fact, $\wis{alg@1} \simeq \wis{commalg}$ by taking as trace maps the identity on every commutative algebra. Further we have a natural commutative diagram of functors
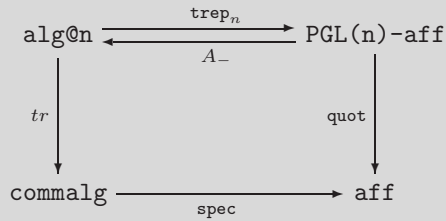
\[
\begin{diagram}
\wis{alg@n} & & \pile{\rTo^{\wis{trep}_n} \\ \lTo_{A_{-}}} & & \wis{PGL(n)-aff} \\
\dTo^{tr} & & & & \dTo^{\wis{quot}} \\
\wis{commalg} & & \rTo_{\wis{spec}} & & \wis{aff}
\end{diagram}
\]
where the bottom map is the equivalence between affine algebras and affine schemes and the top map is the correspondence between Cayley-Hamilton algebras and affine $PGL_n$-schemes, which is {\em not} an equivalence of categories.
\end{boxitpara}

\section{Smooth orders}

To finish this talk let us motivate and define the notion of a {\em Smooth order} properly. Among the many characterizations of commutative regular algebras is the following due to A. Grothendieck.

\begin{theorem} A commutative $\C$-algebra $A$ is regular if and only if it satisfies the following lifting property : if $(B,I)$ is a test-object such that $B$ is a commutative algebra and $I$ is a nilpotent ideal of $B$, then for any algebra map $\phi$, there exists a lifted algebra morphism $\tilde{\phi}$
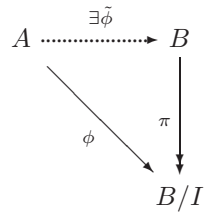
\[
\begin{diagram}
A & \rDotsto^{\exists \tilde{\phi}} & B \\
& \rdTo_{\phi} & \dOnto^{\pi} \\
& & B/I \end{diagram}
\]
 making the diagram commutative.
 \end{theorem}
 
As the category $\wis{commalg}$ of all commutative $\C$-algebras is just $\wis{alg@1}$  it makes sense to define Smooth Cayley-Hamilton algebras by the same lifting property. This was done first by W. Schelter \cite{Schelter} in the category of all algebras satisfying all polynomial identities of $n \times n$ matrices and later by C. Procesi \cite{ProcesiCH} in $\wis{alg@n}$.

\begin{definition} A {\em Smooth Cayley-Hamilton algebra} $A$ is an object in $\wis{alg@n}$ satisfying the following lifting property. If $(B,I)$ is a test-object in $\wis{alg@n}$, that is, $B$ is an object in $\wis{alg@n}$, $I$ is a nilpotent ideal in $B$ such that $B/I$ is an object in $\wis{alg@n}$ and such that the natural map $B \rOnto^{\pi} B/I$ is trace preserving, then every trace preserving algebra map $\phi$ has a lift $\tilde{\phi}$
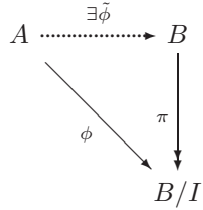
 \[
\begin{diagram}
A & \rDotsto^{\exists \tilde{\phi}} & B \\
& \rdTo_{\phi} & \dOnto^{\pi} \\
& & B/I \end{diagram}
\]
 making the diagram commutative. If $A$ is in addition an order, we say that $A$ is a {\em Smooth order}.
\end{definition}

Next talk we will give a large class of Smooth orders but again it should be stressed that there is no connection between this notion of non-commutative smoothness and the more homological notion of Regular orders (except in dimension one when all notions coincide). 

Still, in the context of $PGL_n$-equivariant affine geometry this notion of non-commutative smoothness is quite natural as illustrated by the following result due to C. Procesi \cite{ProcesiCH}.

\begin{theorem} An object $A$ in $\wis{alg@n}$ is Smooth if and only if the corresponding affine $PGL_n$-scheme $\wis{trep}_n~A$ is smooth (and hence, in particular, reduced).
\end{theorem}

\begin{proof} (One implication) Assume $A$ is Smooth, then to prove that $\wis{trep}_n~A$ is smooth we have to prove that $\C[\wis{trep}_n~A]$ satisfies Grothendieck's  lifting property. So let $(B,I)$ be a test-object in $\wis{commalg}$ and take an algebra morphism $\phi~:~\C[\wis{trep}_n~A] \rTo B/I$. Consider the following diagram
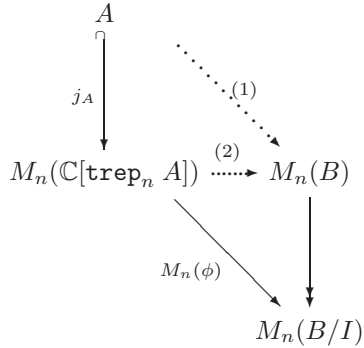
\[
\begin{diagram}
A & & \\
\dInto^{j_A} & \rdDotsto^{(1)} & \\
M_n(\C[\wis{trep}_n~A]) & \rDotsto^{(2)} & M_n(B) \\
& \rdTo_{M_n(\phi)} & \dOnto \\
& & M_n(B/I)
\end{diagram}
\]
the morphism $(1)$ follows from Smoothness of $A$ applied to the morphism $M_n(\phi) \circ j_A$. From the universal property of the map $j_A$ it follows that there is a morphism $(2)$ which is of the form $M_n(\psi)$ for some algebra morphism $\psi~:~\C[\wis{trep}_n~A] \rTo B$. This $\psi$ is the required lift.
\end{proof}

\begin{mexample} Trace rings $\mathbb{T}_{m,n}$ are the free algebras generated by $m$ elements in $\wis{alg@n}$ and as such trivially satisfy the lifting property so are Smooth orders. Alternatively, because
\[
\wis{trep}_n~\mathbb{T}_{m,n} \simeq M_n(\C) \oplus \hdots \oplus M_n(\C) = \C^{mn^2} \]
is a smooth $PGL_n$-variety, $\mathbb{T}_{m,n}$ is Smooth by the previous result.
\end{mexample}

\begin{example} Any commutative algebra $C$ can be viewed as an element of $\wis{alg@n}$ via the diagonal embedding $C \rInto M_n(C)$. However, if $C$ is a regular commutative algebra it is {\em not} true that $C$ is Smooth in $\wis{alg@n}$. For example, take $C = \C[x_1,\hdots,x_d]$ and consider the $4$-dimensional non-commutative local algebra
\[
B = \frac{\C \langle x,y \rangle}{(x^2,y^2,xy+yx)} = \C \oplus \C x \oplus \C y \oplus \C xy \]
with the obvious trace map so that $B \in \wis{alg@2}$. $B$ has a nilpotent ideal $I = B(xy-yx)$ such that the quotient $B/I$ is a $3$-dimensional commutative algebra. Consider the algebra map
\[
\C[x_1,\hdots,x_d] \rTo^{\phi} \frac{B}{I} \qquad \text{defined by} \qquad x_1 \mapsto x \quad x_2 \mapsto y \quad \text{and} \quad x_i \mapsto 0 \quad \text{for $i \geq 3$} \]
This map has no lift as for any potential lifted morphism $\tilde{\phi}$ we have
\[
[ \tilde{\phi}(x),\tilde{\phi}(y) ] \not= 0 \]
whence $\C[x_1,\hdots,x_d]$ is not Smooth in $\wis{alg@2}$.
\end{example}

\begin{mexample} Consider again the quantum plane at $-1$
\[
A = \frac{\C \langle x,y \rangle}{(xy+yx)} \]
then we have seen that $\wis{trep}_2~A = \V(bf+ce) \subset \Af^6$ has a unique isolated singularity at the origin. Hence, $A$ is not a Smooth order.
\end{mexample}

\begin{boxitpara}{box 0.85 setgray fill}
{\bf jotter  : }

Under the correspondence between $\wis{alg@n}$ and $\wis{PGL(n)-aff}$, Smooth Cayley-Hamilton algebras correspond to smooth $PGL_n$-varieties.
\end{boxitpara}

\chapter{non-commutative geometry}

Last time we introduced $\wis{alg@n}$ as a level $n$ generalization of $\wis{commalg}$, the variety of all commutative algebras. Today we will associate to any $A \in \wis{alg@n}$ a {\em non-commutative variety} $\wis{max}~A$ and argue that this gives a non-commutative manifold if $A$ is a Smooth order. In particular we will show that for fixed $n$ and central dimension $d$ there are a finite number of \'etale types of such orders. This fact is the non-commutative analogon of the fact that every manifold is locally diffeomorphic to affine space or, in ringtheory terms, that the $\mathfrak{m}$-adic completion of a regular algebra $C$ of dimension $d$ has just one \'etale type : $\hat{C}_{\mathfrak{m}} \simeq \C [[ x_1,\hdots,x_d ]]$.

\section{Why non-commutative geometry? }

\begin{boxitpara}{box 0.85 setgray fill}
{\bf jotter  : }

There is one new feature that non-commutative geometry has to offer compared to commutative geometry : distinct points can lie infinitesimally close to each other. As desingularization is the process of separating bad tangents, this fact should be useful somehow in our project.
\end{boxitpara}

Recall that if $X$ is an affine commutative variety with coordinate ring $R$, then to each point $P \in X$ corresponds a maximal ideal $\mathfrak{m}_P \triangleleft R$ and a one-dimensional simple representation
\[
S_P = \frac{R}{\mathfrak{m}_P} \]
A basic tool in the study of Hilbert schemes is that  finite closed subschemes of $X$ can be decomposed according to their support. In algebraic terms this means that there are no extensions between different points, that if $P \not= Q$ then
\[
Ext^1_{R}(S_P,S_Q) = 0 \qquad \text{whereas} \qquad Ext^1_R(S_P,S_P) = T_P~X \]
In more plastic lingo : all infinitesimal information of $X$ near $P$ is contained in the self-extensions of $S_P$ and distinct points do not contribute. This is no longer the case for non-commutative algebras.

\begin{mexample} Take the path algebra $A$ of the {\em quiver} $\xymatrix{\vtx{} & \vtx{} \ar[l]}$, that is
\[
A \simeq \begin{bmatrix} \C & \C \\ 0 & \C \end{bmatrix}
\]
Then $A$ has two maximal ideals and two corresponding one-dimensional simple representations
\[
S_1 = \begin{bmatrix} \C \\ 0 \end{bmatrix} = \begin{bmatrix} \C & \C \\ 0 & \C \end{bmatrix}/ \begin{bmatrix} 0 & \C \\ 0 & \C \end{bmatrix} \qquad \text{and} \qquad S_2 = \begin{bmatrix} 0 \\ \C \end{bmatrix} = \begin{bmatrix} \C & \C \\ 0 & \C \end{bmatrix} / \begin{bmatrix} \C & \C \\ 0 & 0 \end{bmatrix} \]
Then, there is a non-split exact sequence with middle term the second column of $A$
\[
0 \rTo S_1= \begin{bmatrix} \C \\ 0 \end{bmatrix} \rTo M = \begin{bmatrix} \C \\ \C \end{bmatrix} \rTo S_2 = \begin{bmatrix} 0 \\ \C \end{bmatrix} \rTo 0 \]
Whence $Ext^1_{A}(S_2,S_1) \not= 0$ whereas $Ext^1_A(S_1,S_2) = 0$. It is no accident that these two facts are encoded into the quiver.
\end{mexample}

\begin{definition} For $A$ an algebra in $\wis{alg@n}$, define its {\em maximal ideal spectrum} $\wis{max}~A$ to be
 the set of all maximal twosided ideals $M$ of $A$ equipped with the {\em non-commutative Zariski topology}, that is, a typical open set of $\wis{max}~A$ is of the form
\[
\mathbb{X}(I) = \{ M \in \wis{max}~A~|~I \not\subset M \} \]
Recall that for every $M \in \wis{max}~A$ the quotient 
\[
\frac{A}{M} \simeq M_k(\C) \qquad \text{for some $k \leq n$} \]
that is, $M$ determines a unique $k$-dimensional simple representation $S_M$ of $A$. 
\end{definition}

As every maximal ideal $M$ of $A$ intersects the center $R$ in a maximal ideal $\mathfrak{m}_P = M \cap R$ we get, in the case of an $R$-order $A$ a continuous map
\[
\wis{max}~A \rTo^c X \quad \text{defined by} \quad M \mapsto P \quad \text{where $M \cap R = \mathfrak{m}_P$} \]
Ringtheorists have studied the fibers $c^{-1}(P)$ of this map in the seventies and eighties in connection with localization theory. The oldest description is the {\em Bergman-Small} theorem, see for example \cite{BergmanSmall} 

\begin{theorem}[Bergman-Small] If $c^{-1}(P) = \{ M_1,\hdots,M_k \}$ then there are natural numbers $e_i \in \N_+$ such that
\[
n = \sum_{i=1}^k e_id_i \qquad \text{where $d_i = dim_{\C}~S_{M_i}$} \]
In particular, $c^{-1}(P)$ is finite for all $P$.
\end{theorem}

Here is a modern proof of this result based on the results of the previous lecture. Because $X$ is the algebraic quotient $\wis{trep}_n~A // GL_n$, points of $X$ correspond to {\em closed} $GL_n$-orbits in $\wis{rep}_n~A$. By a result of M. Artin \cite{Artin69} closed orbits are precisely the isomorphism classes of {\em semi-simple} $n$-dimensional representations, and therefore we denote the quotient variety
\[
X = \wis{trep}_n~A // GL_n = \wis{tiss}_n~A \]
So, a point $P$ determines a semi-simple $n$-dimensional $A$-representation
\[
M_P = S_1^{\oplus e_1} \oplus \hdots \oplus S_k^{\oplus e_k} \]
with the $S_i$ the distinct simple components, say of dimension $d_i = dim_{\C}~S_i$ and occurring in $M_P$ with multiplicity $e_i \geq 1$. This gives $n = \sum e_id_i$ and clearly the annihilator of $S_i$ is a maximal ideal $M_i$ of $A$ lying over $\mathfrak{m}_P$. 

Another interpretation of $c^{-1}(P)$ follows from the work of A. V. Jategaonkar and B. M\"uller. Define a {\em link diagram} on the points of $\wis{max}~A$ by the rule
\[
M \rightsquigarrow M' \qquad \Leftrightarrow \qquad Ext^1_A(S_M,S_{M'}) \not= 0 \]
In fancier language, $M \rightsquigarrow M'$ if and only if $M$ and $M'$ lie infinitesimally close together in $\wis{max}~A$. In fact, the definition of the link diagram in \cite[Chp. 5]{Jategaonkar} or \cite[Chp. 11]{GoodearlWarfield} is slightly different but amounts to the same thing.

\begin{theorem}[Jategaonkar-M\"uller] The connected components of the link diagram on $\wis{max}~A$ are all finite and are in one-to-one correspondence with $P \in X$. That is, if
\[
\{ M_1,\hdots,M_k \} = c^{-1}(P) \subset \wis{max}~A \]
then this set is a connected component of the link diagram.
\end{theorem}

\begin{boxitpara}{box 0.85 setgray fill}
{\bf jotter  : }

In $\wis{max}~A$ there is a Zariski open set of {\em Azumaya points}, that is those $M \in \wis{max}~A$ such that $A/M \simeq M_n(\C)$. It follows that each of these maximal ideals is a singleton connected component of the link diagram. So on this open set there is a one-to-one correspondence between points of $X$ and maximal ideals of $A$ so we can say that $\wis{max}~A$ and $X$ are {\em birational}. However, over the ramification locus there may be several maximal ideals of $A$ lying over the same central maximal ideal and these points should be thought of as lying infinitesimally close to each other. 
\end{boxitpara}
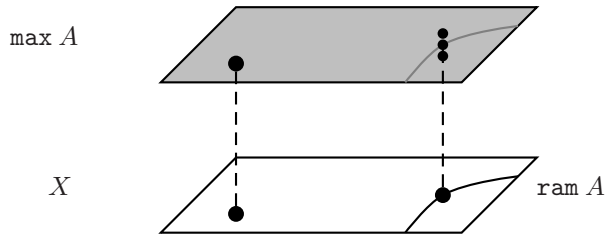
\[
\begin{pspicture}(-2,0)(7.5,3)
\psset{xunit=.5cm}
\psset{yunit=.5cm}
\pspolygon[fillstyle=solid,fillcolor=lightgray](0,4)(2,6)(10,6)(8,4)
\pspolygon[fillstyle=solid,fillcolor=white](0,0)(2,2)(10,2)(8,0)
\pscurve[linecolor=gray](6.5,4)(7.5,5)(9.5,5.5)
\pscurve(6.5,0)(7.5,1)(9.5,1.5)
\put(5,.5){$\wis{ram}~A$}
\put(-1.5,.5){$X$}
\put(-2,2.5){$\wis{max}~A$}
\cnode*(2,4.5){3pt}{A}
\cnode*(2,0.5){3pt}{B}
\cnode*(7.5,5){2pt}{C}
\cnode*(7.5,5.3){2pt}{C1}\cnode*(7.5,4.7){2pt}{C2}
\cnode*(7.5,1){3pt}{D}
\ncline[linestyle=dashed]{A}{B}\ncline[linestyle=dashed]{C}{D}
\end{pspicture}
\]
One might hope that the cluster of infinitesimally points of $\wis{max}~A$ lying over a central singularity $P \in X$ can be used to separate tangent information in $P$ rather than having to resort to the blowing-up process to achieve this.

\section{What non-commutative geometry?}

As an $R$-order $A$ in a central simple $K$-algebra $\Sigma$ of dimension $n^2$ is a finite $R$-module, we can associate to $A$ the sheaf $\Oscr_A$ of non-commutative $\Oscr_X$-algebras using central localization. That is, the section over a basic affine open piece $\mathbb{X}(f) \subset X$ are
\[
\Gamma(\mathbb{X}(f),\Oscr_A) = A_f = A \otimes_R R_f \]
which is readily checked to be a sheaf with global sections $\Gamma(X,\Oscr_A) = A$. As we will investigate Smooth orders via their (central) \'etale structure, that is information about $\hat{A}_{\mathfrak{m}_P}$, we will only need the structure sheaf $\Oscr_A$ over $X$.

In the '70-ties F. Van Oystaeyen \cite{FVO444} and A. Verschoren \cite{LNM887} introduced genuine non-commutative structure sheaves associated to an $R$-order $A$. It is not my intention to promote nostalgia here but perhaps these non-commutative structure sheaves $\Oscr^{nc}_A$ on $\wis{max}~A$ deserve renewed investigation.

\begin{definition}
$\Oscr_A^{nc}$ is defined by taking as the sections over the typical open set $\mathbb{X}(I)$ (for $I$ a twosided ideal of $A$) in $\wis{max}~A$
\[
\Gamma(\mathbb{X}(I),\Oscr_A^{nc}) = \{ \delta \in \Sigma~|~\exists l \in \N~:~I^l \delta \subset A~\} \]
By \cite{FVO444} this defines a sheaf of non-commutative algebras over $\wis{max}~A$ with global sections $\Gamma(\wis{max}~A,\Oscr^{nc}_A) = A$. The stalk of this sheaf at a point $M \in \wis{max}~A$ is the {\em symmetric localization}
\[
\Oscr_{A,M}^{nc} = Q_{A-M}(A) = \{ \delta \in \Sigma~|~I \delta \subset A~\text{for some ideal $I \not\subset P$}~\} \]
\end{definition}

Often, these stalks have no pleasant properties but in some examples, these non-commutative stalks are nicer than those of the central structure sheaf.

\begin{mexample} \label{hereditary} Let $X = \Af^1$, that is, $R = \C[x]$ and consider the order
\[
A = \begin{bmatrix} R & R \\ \mathfrak{m} & R \end{bmatrix} \]
where $\mathfrak{m} = (x) \triangleleft R$. $A$ is an hereditary order so is both a Regular order and a Smooth order. The ramification locus of $A$ is $P_0 = \V(x)$ so over any $P_0 \not= P \in \Af^1$ there is a unique maximal ideal of $A$ lying over $\mathfrak{m}_P$ and the corresponding quotient is $M_2(\C)$. However, over $\mathfrak{m}$ there are two maximal ideals of $A$
\[
M_1 = \begin{bmatrix} \mathfrak{m} & R \\ \mathfrak{m} & R \end{bmatrix} \qquad \text{and} \qquad
M_2 = \begin{bmatrix} R & R \\ \mathfrak{m} & \mathfrak{m} \end{bmatrix} \]
Both $M_1$ and $M_2$ determine a one-dimensional simple representation of $A$, so the Bergman-Small number are $e_1=e_2=1$ and $d_1=d_2=1$. That is, we have the following picture
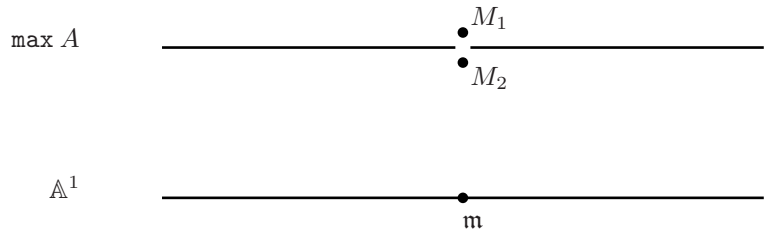
\[
\begin{pspicture}(-2,0)(8,3)
\psline[linewidth=1pt](0,0)(8,0)
\psline[linewidth=1pt](0,2)(3.9,2)
\psline[linewidth=1pt](4.1,2)(8,2)
\cnode*(4,2.2){2pt}{B}
\cnode*(4,1.8){2pt}{C}
\cnode*(4,0){2pt}{D}
\put(4,-.4){$\mathfrak{m}$}
\put(-1.5,0){$\Af^1$}\put(-2,2){$\wis{max}~A$}\put(4.1,2.3){$M_1$}\put(4.1,1.5){$M_2$}
\end{pspicture}
\]

\par \vskip 3mm
There is one non-singleton connected component in the link diagram of $A$ namely
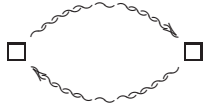
\[
\xymatrix{&& \\ \vtx{}\ar@{~>}@/^4ex/[rr] &&\vtx{}\ar@{~>}@/^4ex/[ll]  \\
&&}
\]
with the vertices corresponding to $\{ M_1,M_2 \}$. The stalk of $\Oscr_{A}$ at the central point $P_0$ is clearly
\[
\Oscr_{A,P_0} = \begin{bmatrix} R_{\mathfrak{m}} & R_{\mathfrak{m}} \\
m_{\mathfrak{m}} & R_{\mathfrak{m}} \end{bmatrix} \]
On the other hand the stalks of the non-commutative structure sheaf $\Oscr_A^{nc}$ in $M_1$ resp. $M_2$ can be computed to be
\[
\Oscr^{nc}_{A,M_1} = \begin{bmatrix} R_{\mathfrak{m}} & R_{\mathfrak{m}} \\ R_{\mathfrak{m}} & R_{\mathfrak{m}} \end{bmatrix} \qquad \text{and} \qquad
\Oscr^{nc}_{A,M_2} = \begin{bmatrix} R_{\mathfrak{m}} & x^{-1}R_{\mathfrak{m}} \\ xR_{\mathfrak{m}} & R_{\mathfrak{m}} \end{bmatrix} \]
and hence both stalks are Azumaya algebras. Observe that we recover the central stalk $\Oscr_{A,P_0}$ as the intersection of these two rings in $M_2(K)$. 

Hence, somewhat surprisingly, the non-commutative structure sheaf of the hereditary non-Azumaya algebra $A$ is a sheaf of Azumaya algebras over $\wis{max}~A$.
\end{mexample}

\section{Marked quiver and Morita settings}

Consider the continuous map for the Zariski topology
\[
\wis{max}~A \rTo^{c} X \]
and let for a central point $P \in X$ the fiber be $\{ M_1,\hdots,M_k \}$ where the $M_i$ are maximal ideals of $A$ with corresponding simple $d_i$-dimensional representation $S_i$. In the previous section we have introduced the {\em Bergman-Small data}, that is 
\[
\alpha = (e_1,\hdots,e_k) \quad \text{and} \quad \beta = (d_1,\hdots,d_k)~\in \N^k_+ \quad \text{satisfying}
\quad \alpha.\beta = \sum_{i=1}^k e_id_i = n \]
(recall that  $e_i$ is the multiplicity of $S_i$ in the semi-simple $n$-dimensional representation corresponding to $P$. Moreover, we have the {\em Jategaonkar-M\"uller data} which is a directed connected graph on the vertices $\{ v_1,\hdots,v_k \}$ (corresponding to the $M_i$) with an arrow
\[
v_i \rightsquigarrow v_j \qquad \text{iff} \qquad Ext^1_A(S_i,S_j) \not= 0 \]
We now want to associate combinatorial objects to this local data. 

To start, introduce a quiver setting $(Q,\alpha)$ where $Q$ is a {\em quiver} (that is, a directed graph) on the vertices $\{ v_1,\hdots,v_k \}$ with the number of arrows from $v_i$ to $v_j$ equal to the dimension of $Ext^1_A(S_i,S_j)$, 
\[
\#~(~v_i \rTo v_j~)~=~dim_{\C}~Ext^1_A(S_i,S_j) \]
and where $\alpha=(e_1,\hdots,e_k)$ is the {\em dimension vector} of the multiplicities $e_i$.

Recall that the representation space $\wis{rep}_{\alpha}~Q$ of a quiver-setting is $\oplus_a M_{e_i \times e_j}(\C)$ where the sum is taken over all arrows $a~:~v_j \rTo v_i$ of $Q$. On this space there is a natural action by the group
\[
GL(\alpha) = GL_{e_1} \times \hdots \times GL_{e_k} \]
by base-change in the vertex-spaces $V_i = \C^{e_i}$ (actually this is an action of $PGL(\alpha)$ which is the quotient of $GL(\alpha)$ by the central subgroup $\C^*(1_{e_1},\hdots,1_{e_k})$). 

The ringtheoretic relevance of the quiver-setting $(Q,\alpha)$ is that
\[
\wis{rep}_{\alpha}~Q \simeq Ext^1_A(M_P,M_P) \qquad \text{as $GL(\alpha)$-modules} \]
where $M_P$ is the semi-simple $n$-dimensional $A$-module corresponding to $P$
\[
M_P = S_1^{\oplus e_1} \oplus \hdots \oplus S_k^{\oplus e_k} \]
and because $GL(\alpha)$ is the automorphism group of $M_P$ there is an induced action on $Ext^1_A(M_P,M_P)$.

Because $M_P$ is $n$-dimensional, an element $\psi \in Ext^1_A(M_P,M_P)$ defines an algebra morphism
\[
A \rTo^{\rho} M_n(\C[\epsilon]) \]
where $\C[\epsilon] = \C[x]/(x^2)$ is the ring of {\em dual numbers}. As we are working in the category $\wis{alg@n}$ we need the stronger assumption that $\rho$ is trace preserving. For this reason we have to consider the $GL(\alpha)$-subspace
\[
tExt^1_A(M_P,M_P) \subset Ext^1_A(M_P,M_P) \]
of {\em trace preserving extensions}. As traces only use blocks on the diagonal (corresponding to loops in $Q$) and as any subspace $M_{e_i}(\C)$ of $\wis{rep}_{\alpha}~Q$ decomposes as a $GL(\alpha)$-module in simple representations
\[
M_{e_i}(\C) = M^0_{e_i}(\C) \oplus \C \]
where $M_{e_i}^0(\C)$ is the subspace of trace zero matrices, we see that
\[
\wis{rep}_{\alpha}~Q^* \simeq tExt^1_A(M_P,M_P) \qquad \text{as $GL(\alpha)$-modules}
\]
where $Q^*$ is a {\em marked quiver} that has the same number of arrows between distinct vertices as $Q$ has, but may have fewer loops and some of these loops may acquire a {\em marking} meaning that their corresponding component in $\wis{rep}_{\alpha}~Q^*$ is $M_{e_i}^0(\C)$ instead of $M_{e_i}(\C)$.

\begin{boxitpara}{box 0.85 setgray fill}
{\bf jotter  : }

Let the local structure of the non-commutative variety $\wis{max}~A$ near the fiber $c^{-1}(P)$ of a point $P \in X$ be determined by the Bergman-Small data
\[
\alpha = (e_1,\hdots,e_k) \qquad \text{and} \qquad \beta = (d_1,\hdots,d_k) \]
and by the Jategoankar-M\"uller data which is encoded in the marked quiver $Q^*$ on $k$-vertices. Then, we associate to $P$ the combinatorial data
\[
\wis{type}(P) = (Q^*,\alpha,\beta) \]
We call $(Q^*,\alpha)$ the {\em marked quiver setting} associated to $A$ in $P \in X$. The dimension vector $\beta = (d_1,\hdots,d_k)$ will be called the {\em Morita setting} associated to $A$ in $P$.
\end{boxitpara}

\begin{mexample} \label{Azumaya} If $A$ is an Azumaya algebra over $R$. then for every maximal ideal $\mathfrak{m}$ corresponding to a point $P \in X$ we have that
\[
A/\mathfrak{m}A = M_n(\C) \]
so there is a unique maximal ideal $M = \mathfrak{m}A$ lying over $\mathfrak{m}$ whence the Jategaonkar-M\"uller data are $\alpha = (1)$ and $\beta = (n)$. If $S_P = R/\mathfrak{m}$ is the simple representation of $R$ we have
\[
Ext^1_A(M_P,M_P) \simeq Ext^1_R(S_P,S_P) = T_P~X \]
and as all the extensions come from the center, the corresponding algebra representations $A \rTo M_n(\C[\epsilon])$ are automatically trace preserving. That is, the marked quiver-setting associated to $A$ in $P$ is
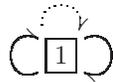
\[
\xymatrix{ & \\ \vtx{1} \ar@(dl,ul) \ar@{.>}@(ul,ur) \ar@(ur,dr) & \\ &} 
\]
where the number of loops is equal to the dimension of the tangent space $T_P~X$ in $P$ at $X$ and the Morita-setting associated to $A$ in $P$ is $(n)$.
\end{mexample}

\begin{mexample} Consider the order of example~\ref{hereditary} which is generated as a $\C$-algebra by the elements
\[
a = \begin{bmatrix} 1 & 0 \\ 0 & 0 \end{bmatrix} \quad b = \begin{bmatrix} 0 & 1 \\ 0 & 0 \end{bmatrix} \quad c = \begin{bmatrix} 0 & 0 \\ x & 0 \end{bmatrix} \quad d = \begin{bmatrix} 0 & 0 \\ 0 & 1 \end{bmatrix} \]
and the $2$-dimensional semi-simple representation $M_{P_0}$ determined by $\mathfrak{m}$ is given by the algebra morphism $A \rTo M_2(\C)$ sending $a$ and $d$ to themselves and $b$ and $c$ to the zero matrix. A calculation shows that 
\[
Ext^1_A(M_{P_0},M_{P_)}) = \wis{rep}_{\alpha}~Q \qquad \text{for} \qquad (Q,\alpha) = \xymatrix{
\vtx{1} \ar@/^/[r]^u & \vtx{1} \ar@/^/[l]^v} \]
and as the correspondence with algebra maps to $M_2(\C[\epsilon])$ is given by
\[
a \mapsto \begin{bmatrix} 1 & 0 \\ 0 & 0 \end{bmatrix} \quad
b \mapsto \begin{bmatrix} 0 & \epsilon v \\ 0 & 0 \end{bmatrix} \quad
c \mapsto \begin{bmatrix} 0 & 0 \\ \epsilon u & 0 \end{bmatrix} \quad
d \mapsto \begin{bmatrix} 0 & 0 \\ 0 & 1 \end{bmatrix} \]
each of these maps is trace preserving so the marked quiver setting is $(Q,\alpha)$ and the Morita-setting is $(1,1)$.
\end{mexample}

\section{Local classification}

\begin{boxitpara}{box 0.85 setgray fill}
{\bf jotter  : }

Because the combinatorial data $\wis{type}(P)=(Q^*,\alpha,\beta)$ encodes the infinitesimal information of the cluster of maximal ideals of $A$ lying over the central point $P \in X$, 
$(\wis{rep}_{\alpha}~Q^*,\beta)$ should be viewed as analogous to the usual tangent space $T_P~X$.

 If $P \in X$ is a singular point, then the tangent space is too large so we have to impose additional relations to describe the variety $X$ in a neighborhood of $P$, but if $P$ is a smooth point we can recover the local structure of $X$ from $T_P~X$.

Here we might expect a similar phenomenon : in general the data $(\wis{rep}_{\alpha}~Q^*,\beta)$ will be too big to describe $\hat{A}_{\mathfrak{m}_P}$ unless $A$ is a Smooth order in $P$ in which case we can recover $\hat{A}_{\mathfrak{m}_P}$.
\end{boxitpara}

We begin by defining some algebras which can be described combinatorially  from $(Q^*,\alpha,\beta)$.

For every arrow $a~:~v_i \rTo v_j$ define a {\em generic rectangular matrix} of size $e_j \times e_i$
\[
X_a = \begin{bmatrix}
x_{11}(a) & \hdots & \hdots & x_{1e_i}(a) \\
\vdots & & & \vdots \\
x_{e_j1}(a) & \hdots & \hdots & x_{e_je_i}(a) \end{bmatrix} \]
(and if $a$ is a marked loop take $x_{e_ie_i}(a) = -x_{11}(a)-x_{22}(a)- \hdots - x_{e_i-1e_i-1}(a)$)
then the coordinate ring $\C[\wis{rep}_{\alpha}~Q^*]$ is the polynomial ring in the entries of all $X_a$. For an oriented path $p$ in the marked quiver $Q^*$ with starting vertex $v_i$ and terminating vertex $v_j$
\[
v_i \rDotsto^p v_j \quad = \quad v_i \rTo^{a_1} v_{i_1} \rTo^{a_2} \hdots \rTo^{a_{l-1}} v_{i_l} \rTo^{a_l} v_j \]
we can form the square $e_j \times e_i$ matrix
\[
X_p = X_{a_l} X_{a_{l-1}} \hdots X_{a_2} X_{a_1} \]
which has all its entries polynomials in $\C[\wis{rep}_{\alpha}~Q^*]$. In particular, if the path is an oriented {\em cycle} $c$ in $Q^*$ starting and ending in $v_i$ then $X_c$ is a square $e_i \times e_i$ matrix and we can take its trace $tr(X_c) \in \C[\wis{rep}_{\alpha}~Q^*]$ which is a polynomial invariant under the action of $GL(\alpha)$ on $\wis{rep}_{\alpha}~Q^*$. 

In fact, it was proved in \cite{LBProcesi} that these {\em traces along oriented cycles} generate the invariant ring
\[
R^{\alpha}_{Q^*} = \C[\wis{rep}_{\alpha}~Q^*]^{GL(\alpha)} \subset \C[\wis{rep}_{\alpha}~Q^*] \]
Next we bring in the Morita-setting $\beta = (d_1,\hdots,d_k)$ and define a block-matrix ring
\[
A^{\alpha,\beta}_{Q^*} = \begin{bmatrix}
M_{d_1 \times d_1}(P_{11}) & \hdots & M_{d_1 \times d_k}(P_{1k}) \\
\vdots & & \vdots \\
M_{d_k \times d_1}(P_{k1}) & \hdots & M_{d_k \times d_k}(P_{kk}) \end{bmatrix} \subset M_n(\C[\wis{rep}_{\alpha}~Q^*]) \]
where $P_{ij}$ is the $R^{\alpha}_{Q^*}$-submodule of $M_{e_j \times e_i}(\C[\wis{rep}_{\alpha}~Q^*])$ generated by all $X_p$ where $p$ is an oriented path in $Q^*$ starting in $v_i$ and ending in $v_k$. 

Observe that for triples $(Q^*,\alpha,\beta_1)$ and $(Q^*,\alpha,\beta_2)$ we have that
\[
A^{\alpha,\beta_1}_{Q^*} \qquad \text{is Morita-equivalent to} \qquad A^{\alpha,\beta_2}_{Q^*} \]
whence the name Morita-setting for $\beta$. 

Before we can state the next result we need the {\em Euler-form} of the underlying quiver $Q$ of $Q^*$ (that is, forgetting the markings of some loops) which is the bilinear form $\chi_Q$ on $\Z^k$ determined by the matrix having as its $(i,j)$-entry $\delta_{ij} - \# \{ a~:~v_i \rTo^a v_j \}$. The statements below can be deduced from those of \cite{LBProcesi}

\begin{theorem} For a triple $(Q^*,\alpha,\beta)$ with $\alpha.\beta = n$ we have 
\begin{enumerate}
\item{$A^{\alpha,\beta}_{Q^*}$ is an $R^{\alpha}_Q$-order in $\wis{alg@n}$ if and only if $\alpha$ is the dimension vector of a simple representation of $Q^*$, that is, for all vertex-dimensions $\delta_i$ we have
\[
\chi_Q(\alpha,\delta_i) \leq 0 \qquad \text{and} \qquad \chi_Q(\delta_i,\alpha) \leq 0 \]
unless $Q^*$ is an oriented cycle of type $\tilde{A}_{k-1}$ then $\alpha$ must be $(1,\hdots,1)$.}
\item{If this condition is satisfied, the dimension of the center $R^{\alpha}_{Q^*}$ is equal to
\[
\wis{dim}~R^{\alpha}_{Q^*} = 1 - \chi_Q(\alpha,\alpha) - \# \{ \text{marked loops in $Q^*$} \} \]
}
\end{enumerate}
\end{theorem}

These combinatorial algebras determine the \'etale local structure of Smooth orders as was proved in \cite{LBLocalStructure}. The principal technical ingredient in the proof is the {\em Luna slice theorem}, see for example \cite{Slodowy} or \cite{Luna}.

\begin{theorem} Let $A$ be a Smooth order over $R$ in $\wis{alg@n}$ and let $P \in X$ with corresponding maximal ideal $\mathfrak{m}$. If the marked quiver setting and the Morita-setting associated to $A$ in $P$ is given by the triple $(Q^*,\alpha,\beta)$, then there is a Zariski open subset $\mathbb{X}(f_i)$ containing $P$ and an \'etale extension $S$ of both $R_{f_i}$ and the algebra $R^{\alpha}_{Q^*}$ such that we have the following diagram
\[
\begin{diagram}
& & A_{f_i} \otimes_{R_{f_i}} S \simeq A^{\alpha,\beta}_{Q^*} \otimes_{R^{\alpha}_{Q^*}} S & & \\
& \ruTo & \uTo & \luTo & \\
A_{f_i} & & S & & A^{\alpha,\beta}_{Q^*} \\
\uTo & \ruTo_{etale} & & \luTo_{etale} & \uTo \\
R_{f_i} & & & & R^{\alpha}_{Q^*} 
\end{diagram}
\]
In particular, we have
\[
\hat{R}_{\mathfrak{m}} \simeq \hat{R}^{\alpha}_{Q^*} \qquad \text{and} \qquad \hat{A}_{\mathfrak{m}} \simeq \hat{A}^{\alpha,\beta}_{Q^*} \]
where the completions at the right hand sides are with respect to the maximal (graded) ideal of $R^{\alpha}_{Q^*}$ corresponding to the zero representation.
\end{theorem}

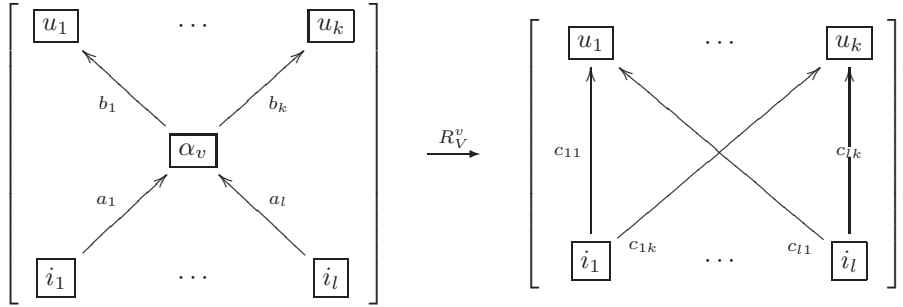
\begin{mexample} From example~\ref{Azumaya} we recall that the triple $(Q^*,\alpha,\beta)$ associated to an Azumaya algebra in a point $P \in X$ is given by

\par \vskip 3mm
\[
\xymatrix{ \vtx{1} \ar@(dl,ul) \ar@{.>}@(ul,ur) \ar@(ur,dr) & } \quad \text{and} \quad \beta = (n)
\]
where the number of arrows is equal to $dim_{\C}~T_P X$. In case $P$ is a smooth point of $X$ this number is equal to $d = \wis{dim}~X$. Observe that $GL(\alpha) = \C^*$ acts trivially on $\wis{rep}_{\alpha}~Q^* = \C^d$ in this case. Therefore we have that
\[
R^{\alpha}_{Q^*} \simeq \C[x_1,\hdots,x_d] \quad \text{and} \quad A^{\alpha,\beta}_{Q^*} = M_n(\C[x_1,\hdots,x_d]) \] 
Because $A$ is a Smooth order in such points we get that 
\[
\hat{A}_{\mathfrak{m_P}} \simeq M_n(\C[[x_1,\hdots,x_d]]) \]
consistent with our \'etale local knowledge of Azumaya algebras.
\end{mexample}

\begin{boxitpara}{box 0.85 setgray fill}
{\bf jotter  : }

Because $\alpha.\beta=n$, the number of vertices of $Q^*$ is bounded by $n$ and as 
\[
d = 1 -\chi_Q(\alpha,\alpha) -\# \{ \text{marked loops} \} \]
the number of arrows and (marked) loops is also bounded. This means that for a particular dimension $d$ of the central variety $X$ there are only a finite number of \'etale local types of Smooth orders in $\wis{alg@n}$. 

This fact might be seen as a non-commutative version of the fact that there is just one \'etale type of a smooth variety in dimension $d$ namely $\C[[x_1,\hdots,x_d]]$. At this moment a similar result for Regular orders seems to be far out of reach.
\end{boxitpara}

\section{A two-person game}

Starting with a marked quiver setting $(Q^*,\alpha)$ we will play a two-person game. Left will be allowed to make one of the reduction steps to be defined below if the condition on Leaving arrows is satisfied, Red on the other hand if the condition on aRRiving arrows is satisfied. Although we will not use combinatorial game theory in any way, it is a very pleasant topic and the interested reader is referred to \cite{ConwayONAG} or \cite{WinningWays}.

The reduction steps below were discovered by R. Bocklandt in his Ph.D. thesis \cite{BocklandtThesis} (see also \cite{Bocklandtpaper}) in which he classifies quiver settings having a regular ring of invariants. These steps were slightly extended in \cite{RBLBVdW} in order to classify central singularities of Smooth orders. All reductions are made locally around a vertex in the marked quiver. There are three types of allowed moves 

{\bf Vertex removal}

Assume we have a marked quiver setting $(Q^*,\alpha)$ and a vertex $v$ such that the local structure of $(Q^*,\alpha)$ near $v$ is indicated by the picture on the left below, that is, inside the vertices we have written the components of the dimension vector and the subscripts of an arrow indicate how many such arrows there are in $Q^*$ between the indicated vertices.
Define the new marked quiver setting $(Q^*_R,\alpha_R)$ obtained by the operation $R^v_V$ which removes the vertex $v$ and composes all arrows through $v$, the dimensions of the other vertices are unchanged :
\[
\left[ ~\vcenter{
\xymatrix@=1cm{
\vtx{u_1}&\cdots &\vtx{u_k}\\
&\vtx{\alpha_v}\ar[ul]^{b_1}\ar[ur]_{b_k}&\\
\vtx{i_1}\ar[ur]^{a_1}&\cdots &\vtx{i_l}\ar[ul]_{a_l}}}
~\right] \quad
\rTo^{R^v_V} \quad
\left[~\vcenter{
\xymatrix@=1cm{
\vtx{u_1}&\cdots &\vtx{u_k}\\
&&\\
\vtx{i_1}\ar[uu]^{c_{11}}\ar[uurr]_<<{c_{1k}}&\cdots &\vtx{i_l}\ar[uu]|{c_{lk}}\ar[uull]^<<{c_{l1}}}}
~\right].
\]
where $c_{ij} = a_ib_j$ (observe that some of the incoming and outgoing vertices may be the
same so that one obtains loops in the corresponding vertex). Left (resp. Right) is allowed to make this reduction step provided the following condition is met
\[
(Left) \quad \chi_Q(\alpha,\epsilon_v) \geq 0 \quad \Leftrightarrow \quad \alpha_v \geq \sum_{j=1}^l a_j i_j \] 
\[
(Right) \quad \chi_Q(\epsilon_v,\alpha) \geq 0\quad \Leftrightarrow \quad \alpha_v \geq \sum_{j=1}^k b_j u_j \]
(observe that if we started off from a marked quiver setting $(Q^*,\alpha)$ coming from an order, then these inequalities must actually be equalities).

{\bf loop removal}

If $v$ is a vertex with vertex-dimension $\alpha_v = 1$ and having $k \geq 1$ loops.
Let $(Q^*_R,\alpha_R)$ be the marked quiver setting obtained by the loop removal operation $R^v_l$
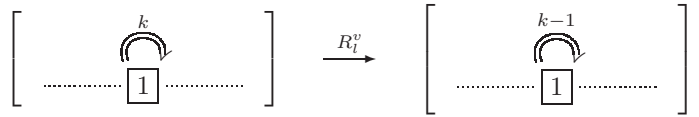
\[
\left[~\vcenter{
\xymatrix@=1cm{
&\vtx{1}\ar@{..}[r]\ar@{..}[l]\ar@(lu,ru)@{=>}^k&}}
~\right]\quad \rTo^{R^v_l} \quad
\left[~\vcenter{
\xymatrix@=1cm{
&\vtx{1}\ar@{..}[r]\ar@{..}[l]\ar@(lu,ru)@{=>}^{k-1}&}}
~\right].\]
removing one loop in $v$ and keeping the same dimension vector. Both Left and Right are allowed to make this reduction step.

{\bf Loop removal}

If the local situation in $v$ is such that there is exactly one (marked) loop in $v$, the dimension vector in $v$ is $k \geq 2$ and there is exactly one arrow Leaving $v$ and this to a vertex with dimension vector $1$, then Left is allowed to make the reduction $R^v_L$ indicated below
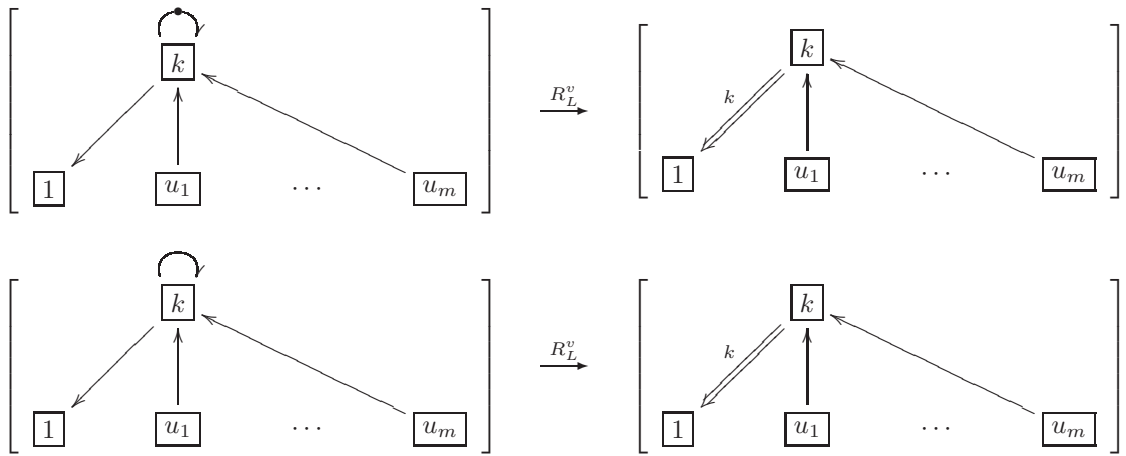
\[
\left[~\vcenter{
\xymatrix@=1cm{
&\vtx{k}\ar[dl]\ar@(lu,ru)|{\bullet}&&\\
\vtx{1}&\vtx{u_1}\ar[u]&\cdots &\vtx{u_m}\ar[ull]}}
~\right]\quad \rTo^{R^v_L} \quad
\left[~\vcenter{
\xymatrix@=1cm{
&\vtx{k}\ar@2[dl]_{k}&&\\
\vtx{1}&\vtx{u_1}\ar[u]&\cdots &\vtx{u_m}\ar[ull]}}
~\right].
\]
\vspace{.5cm}
\[
\left[~\vcenter{
\xymatrix@=1cm{
&\vtx{k}\ar[dl]\ar@(lu,ru)&&\\
\vtx{1}&\vtx{u_1}\ar[u]&\cdots &\vtx{u_m}\ar[ull]}}
~\right]\quad \rTo^{R^v_L} \quad
\left[~\vcenter{
\xymatrix@=1cm{
&\vtx{k}\ar@2[dl]_k&&\\
\vtx{1}&\vtx{u_1}\ar[u]&\cdots &\vtx{u_m}\ar[ull]}}
~\right].
\]

Similarly, if there is one (marked) loop in $v$ and $\alpha_v = k \geq 2$ and there is only one arrow aRRiving at $v$ coming from a vertex of dimension vector $1$, then Right is allowed to make the reduction $R^v_L$
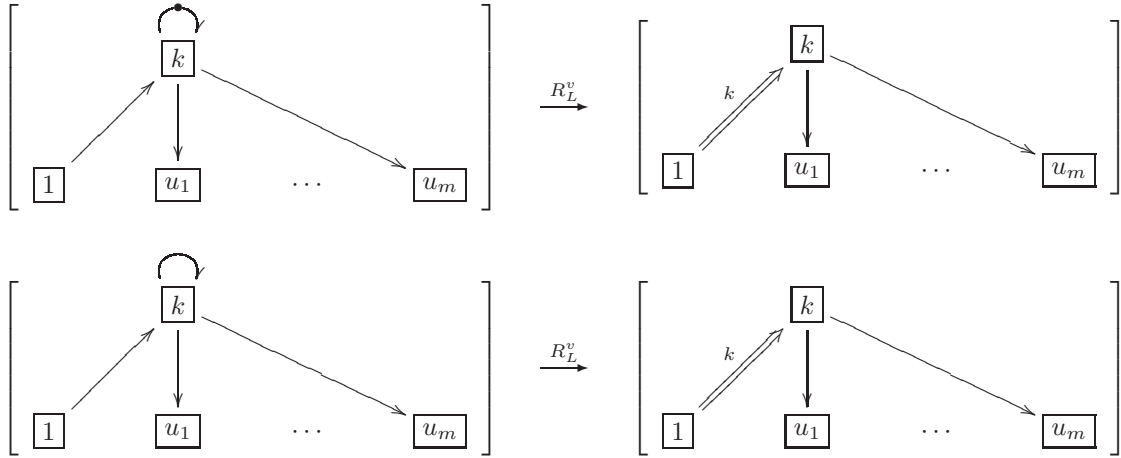
\[
\left[~\vcenter{
\xymatrix@=1cm{
&\vtx{k}\ar[d]\ar[drr]\ar@(lu,ru)|{\bullet}&&\\
\vtx{1}\ar[ur]&\vtx{u_1}&\cdots &\vtx{u_m}}}
~\right]\quad \rTo^{R^v_L} \quad
\left[~\vcenter{
\xymatrix@=1cm{
&\vtx{k}\ar[d]\ar[drr]&&\\
\vtx{1}\ar@2[ur]^k&\vtx{u_1}&\cdots &\vtx{u_m}}}
~\right] \]
\vspace{.5cm}
\[
\left[~\vcenter{
\xymatrix@=1cm{
&\vtx{k}\ar[d]\ar[drr]\ar@(lu,ru)&&\\
\vtx{1}\ar[ur]&\vtx{u_1}&\cdots &\vtx{u_m}}}
~\right]\quad \rTo^{R^v_L} \quad
\left[~\vcenter{
\xymatrix@=1cm{
&\vtx{k}\ar[d]\ar[drr]&&\\
\vtx{1}\ar@2[ur]^k&\vtx{u_1}&\cdots &\vtx{u_m}}}
~\right]
\]

In accordance with combinatorial game theory we call a marked quiver setting $(Q^*,\alpha)$ a {\em zero setting} if neither Left nor Right has a legal reduction step. The relevance of this game on marked quiver settings is that if
\[
(Q^*_1,\alpha_1) \rightsquigarrow (Q^*_2,\alpha_2) \]
is a sequence of legal moves (both Left and Right are allowed to pass), then 
\[
R^{\alpha_1}_{Q^*_1} \simeq R^{\alpha_2}_{Q^*_2}[y_1,\hdots,y_z] \]
where $z$ is the sum of all loops removed in $R^v_l$ reductions plus the sum of $\alpha_v$ for each reduction step $R^v_L$ involving a genuine loop and the sum of $\alpha_v - 1$ for each reduction step $R^v_L$ involving a marked loop. That is, marked quiver settings which below to the same game tree have smooth equivalent invariant rings.

In general games, a position can reduce to several zero-positions depending on the chosen moves. For this reason the next result, proved in \cite{RBLBVdW} is somewhat surprising

\begin{theorem} Let $(Q^*,\alpha)$ be a marked quiver setting, then there is a unique zero-setting
$(Q^*_0,\alpha_0)$ for which there exists a reduction procedure
\[
(Q^*,\alpha) \rightsquigarrow (Q^*_0,\alpha_0) \]
We will denote this unique zero-setting by $Z(Q^*,\alpha)$.
\end{theorem}

\begin{boxitpara}{box 0.85 setgray fill}
{\bf jotter  : }

Therefore it is sufficient to classify the zero-positions if we want to characterize all central singularities of a Smooth order in a given central dimension $d$.
\end{boxitpara}

\section{Central singularities}

Let $A$ be a Smooth $R$-order in $\wis{alg@n}$ and $P$ a point in the central variety $X$ with corresponding maximal ideal $\mathfrak{m} \triangleleft R$.  We now want to classify the types of singularities of $X$ in $P$, that is to classify $\hat{R}_{\mathfrak{m}}$. 

To start, can we decide when $P$ is a smooth point of $X$ ?  In the case that $A$ is an Azumaya algebra in $P$, we know already that $A$ can only be a Smooth if $R$ is regular in $P$. Moreover we have seen for $A$ a Regular reflexive Azumaya algebra that the non-Azumaya points in $X$ are precisely the singularities of $X$. 

For Smooth orders the situation is more delicate but as mentioned before we have a complete solution in terms of the two-person game by a slight adaptation of Bocklandt's main result \cite{Bocklandtpaper}.

\begin{theorem} If $A$ is a Smooth $R$-order and $(Q^*,\alpha,\beta)$ is the combinatorial data associated to $A$ in $P \in X$. Then, $P$ is a smooth point of $X$ if and only if the unique associated zero-setting
\[
Z(Q^*,\alpha) \in \{~\xymatrix{\vtx{1}} \qquad \xymatrix{\vtx{k} \ar@(ul,ur)}  \qquad  \xymatrix{\vtx{k} \ar@(ul,ur)|{\bullet}}  \quad~\qquad  \xymatrix{\vtx{2} \ar@(dl,ul) \ar@(dr,ur)} \qquad~\quad~\quad  \xymatrix{\vtx{2} \ar@(dl,ul) \ar@(dr,ur)|{\bullet}}~\quad \xymatrix{\vtx{2} \ar@(dl,ul)|{\bullet} \ar@(dr,ur)|{\bullet}} \qquad~\}
\]

The Azumaya points are such that $Z(Q^*,\alpha) = \xymatrix{\vtx{1}}$ hence the singular locus of $X$ is contained in the ramification locus $\wis{ram}~A$ but may be strictly smaller.
\end{theorem}

To classify the central singularities of Smooth orders we may reduce to zero-settings $(Q^*,\alpha) = Z(Q^*,\alpha)$. For such a setting we have for all vertices $v_i$ the inequalities
\[
\chi_Q(\alpha,\delta_i) < 0 \qquad \text{and} \qquad \chi_Q(\delta_i,\alpha) < 0 \]
and the dimension of the central variety can be computed from the Euler-form $\chi_Q$. This gives us an estimate of $d = \wis{dim}~X$ which is very efficient to classify the singularities in low dimensions.

\begin{theorem} \label{counting} Let $(Q^*,\alpha) = Z(Q^*,\alpha)$ be a zero-setting on $k \geq 2$ vertices. Then,
\[
\wis{dim}~X  \geq 1 + \sum_{\xymatrix@=1cm{ \vtx{a} }}^{a \geq 1} a + 
\sum_{\xymatrix@=1cm{ \vtx{a}\ar@(ul,dl)|{\bullet} }}^{a > 1}(2a-1) +
 \sum_{\xymatrix@=1cm{ \vtx{a}\ar@(ul,dl)}}^{a > 1}(2a) + \sum_{\xymatrix@=1cm{ \vtx{a}\ar@(ul,dl)|{\bullet}\ar@(ur,dr)|{\bullet}}}^{a > 1} (a^2+a-2) + \]
\[
\sum_{\xymatrix@=1cm{ \vtx{a}\ar@(ul,dl)|{\bullet}\ar@(ur,dr)}}^{a > 1} (a^2+a-1) +
\sum_{\xymatrix@=1cm{ \vtx{a}\ar@(ul,dl)\ar@(ur,dr)}}^{a > 1} (a^2+a) + \hdots +
\sum_{\xymatrix@=1cm{ \vtx{a}\ar@(ul,dl)|{\bullet}_{k}\ar@(ur,dr)^{l}}}^{a > 1} ((k+l-1)a^2+a-k) + \hdots
\]
In this sum the contribution of a vertex $v$ with $\alpha_v = a$ is determined by the number of
(marked) loops in $v$. By the reduction steps (marked) loops only occur at vertices where
$\alpha_v > 1$.
\end{theorem}

Let us illustrate this result by classifying the central singularities in low dimensions

\begin{mexample}[dimension $2$] When $\wis{dim}~X = 2$ no zero-position on at least two vertices satisfies the inequality of theorem~\ref{counting}, so the only zero-position possible to be obtained from a marked quiver-setting $(Q^*,\alpha)$ in dimension two is
\[
Z(Q^*,\alpha) = \xymatrix{\vtx{1}} \]
and therefore the central two-dimensional variety $X$ of a Smooth order is smooth.
\end{mexample}

\begin{mexample}[dimension $3$] If $(Q^*,\alpha)$ is a zero-setting for dimension $\leq 3$ then $Q^*$ can have at most two vertices. If there is just one vertex it must have dimension $1$ (reducing again to $\xymatrix{\vtx{1}}$ whence smooth) or must be
\[
Z(Q^*,\alpha) = \qquad \xymatrix{\vtx{2} \ar@(ul,dl)|{\bullet} \ar@(ur,dr)|{\bullet}} \]
which is again a smooth setting. If there are two vertices both must have dimension $1$ and both must have at least two incoming and two outgoing arrows (for otherwise we could perform an additional vertex-removal reduction). As there are no loops possible in these vertices for zero-settings, it follows from the formula $d = 1 - \chi_Q(\alpha,\alpha)$ that the only possibility is
\[
Z(Q^*,\alpha) = \xymatrix{\vtx{1} \ar@/^2ex/[rr]_a \ar@/^4ex/[rr]^b & & \vtx{1} \ar@/^2ex/[ll]_c \ar@/^4ex/[ll]^d} \]
The ring of polynomial invariants $R^{\alpha}_{Q^*}$ is generated by traces along oriented cycles in $Q^*$ so in this case it is generated by the invariants
\[
x = ac, \quad y = ad, \quad u = bc \quad \text{and} \quad v = bd \]
and there is one relation between these generators, so
\[
R^{\alpha}_{Q^*} \simeq \frac{\C[x,y,u,v]}{(xy-uv)} \]
Therefore, the only \'etale type of central singularity in dimension three is the {\em conifold singularity}.
\end{mexample}

\begin{mexample}[dimension $4$] If $(Q^*,\alpha)$ is a zero-setting for dimension $4$ then $Q^*$ can have at most three vertices. If there is just one, its dimension must be $1$ (smooth setting) or $2$ in which case the only new type is
\[
Z(Q^*,\alpha) = \qquad \xymatrix{\vtx{2} \ar@(ul,dl) \ar@(ur,dr)|{\bullet}} \]
which is again a smooth setting.

If there are two vertices, both must have dimension $1$ and have at least two incoming and outgoing arrows as in the previous example. The only new type that occurs is

\par \vskip 3mm
\[
Z(Q^*,\alpha) = \xymatrix{ \vtx{1} \ar@/^/[rr] \ar@/^3ex/[rr] & & \vtx{1} \ar@/^/[ll] \ar@/^2ex/[ll] \ar@/^3ex/[ll]} \]

for which one calculates as before the ring of invariants to be
\[
R^{\alpha}_{Q^*} = \frac{\C[a,b,c,d,e,f]}{(ae-bd,af-cd,bf-ce)} \]
If there are three vertices all must have dimension $1$ and each vertex must have at least two incoming and two outgoing vertices. There are just two such possibilities in dimension $4$
\[
Z(Q^*,\alpha) \in \{ \quad 
\xymatrix{\vtx{1}\ar@/^/[rr]\ar@/^/[rd]&&\vtx{1}\ar@/^/[ll]\ar@/^/[ld]\\
&\vtx{1}\ar@/^/[ru]\ar@/^/[lu]&}  \qquad 
\xymatrix{\vtx{1}\ar@2@/^/[rr]&&\vtx{1}\ar@2@/^/[ld]\\
&\vtx{1}\ar@2@/^/[lu]&} \quad \}
\]
The corresponding rings of polynomial invariants are
\[
R^{\alpha}_{Q^*} = \frac{\C[x_1,x_2,x_3,x_4,x_5]}{(x_4x_5-x_1x_2x_3)} \qquad \text{resp.} \qquad
R^{\alpha}_{Q^*} = \frac{\C[x_1,x_2,x_3,x_4,y_1,y_2,y_3,y_4]}{R_2} \]
where $R_2$ is the ideal generated by all $2 \times 2$ minors of the matrix
\[
\begin{bmatrix} 
x_1 & x_2 & x_3 & x_4 \\
y_1 & y_2 & y_3 & y_4 \end{bmatrix}
\]
\end{mexample}

In \cite{RBLBVdW} it was proved that there are exactly ten types of Smooth order central singularities in dimension $d=5$ and $53$ in dimension $d=6$. The strategy to prove such a result  is as follows. 

First one makes a full list of all zero-settings $(Q^*,\alpha) = Z(Q^*,\alpha)$ such that $d = 1 -\chi_Q(\alpha,\alpha) - \#$ marked loops, using theorem~\ref{counting}. 

Next, one has to weed out zero-settings having isomorphic rings of polynomial invariants (or rather, having the same $\mathfrak{m}$-adic completion where $\mathfrak{m} \triangleleft R^{\alpha}_{Q^*}$ is the unique graded maximal ideal generated by all generators). There are two invariants to separate two rings of invariants. 

One is the sequence of numbers
\[
dim_{\C}~\frac{\mathfrak{m}^n}{\mathfrak{m}^{n+1}} \]
which can sometimes be computed easily (for example if all dimension vector components are equal to $1$). 

The other invariant is what we call the {\em fingerprint} of the singularity. In most cases, there will be other types of singularities (necessarily also of Smooth order type) in the variety corresponding to $R^{\alpha}_{Q^*}$ and the methods of \cite{LBLocalStructure} allow us to determine their associated marked quiver settings as well as the dimensions of these strata. 

In most cases these two methods allow to separate the different types of singularities. In the few remaining cases it is then easy to write down an explicit isomorphism. We refer to (the published version of) \cite{RBLBVdW} for the full classification of these singularities in dimension $5$ and $6$.

\begin{boxitpara}{box 0.85 setgray fill}
{\bf jotter  : }

In low dimensions there is a full classification of all central singularities $\hat{R}_{\mathfrak{m}}$ of a Smooth order in $\wis{alg@n}$. However, at this moment no such classification exists for $\hat{A}_{\mathfrak{m}}$. That is, under the game rules it is not clear what structural results of the orders $A^{\alpha}_{Q^*}$ are preserved.
\end{boxitpara}

\section{Isolated singularities}

In the classification of central singularities of Smooth orders, isolated singularities stand out as the fingerprinting method to separate them clearly fails. Fortunately, we do have by \cite{RBLBSS} a complete classification of these (in all dimensions).

\begin{theorem} Let $A$ be a Smooth order over $R$ and let $(Q^*,\alpha,\beta)$ be the combinatorial data associated to a $A$ in a point $P \in X$. Then, $P$ is an isolated singularity if and only if $Z(Q^*,\alpha) = T(k_1,\hdots,k_l)$ where
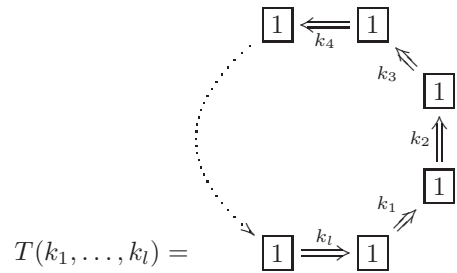
\[
T(k_1,\hdots,k_l) = \xy 0;/r.15pc/:
\POS (0,0) *+{\vtx{1}} ="a",
(20,0) *+{\vtx{1}} ="b",
(34,14) *+{\vtx{1}} ="c",
(34,34) *+{\vtx{1}} ="d",
(20,48) *+{\vtx{1}} ="e",
(0,48) *+{\vtx{1}} ="f"
\POS"a" \ar@{=>}^{k_l} "b"
\POS"b" \ar@{=>}^{k_1} "c"
\POS"c" \ar@{=>}^{k_2} "d"
\POS"d" \ar@{=>}^{k_3} "e"
\POS"e" \ar@{=>}^{k_4} "f"
\POS"f" \ar@/_7ex/@{.>} "a"
\endxy 
\]
with $d = \wis{dim}~X = \sum_i k_i - l +1$. 

Moreover, two such singularities, corresponding to $T(k_1,\hdots,k_l)$ and $T(k'_1,\hdots,k'_{l'})$,  are isomorphic if and only if 
\[
l = l' \qquad \text{and} \qquad  k'_i = k_{\sigma(i)} \]
for some permutation $\sigma \in S_l$.
\end{theorem}

The results we outlined in this talk are good as well as bad news. 

\begin{boxitpara}{box 0.85 setgray fill}
{\bf jotter  : }

On the positive side we have very precise information on the types of singularities which can occur in the central variety of a Smooth order (certainly in low dimensions) in sharp contrast to the case of Regular orders. 

However, because of the scarcity of such types most interesting quotient singularities $\C^d/G$ will {\em not} have a Smooth order over their coordinate ring $R = \C[\C^d/G]$. 
\end{boxitpara}

So, after all this hard work we seem to have come to a dead end with respect to the desingularization problem as there are no Smooth orders with center $C[\C^d/G]$. Fortunately, we have one remaining trick available : to bring in a {\em stability structure}.

\chapter{non-commutative desingularizations}

In the first talk I claimed that in order to find good desingularizations of quotient singularities $\C^d/G$ we had to find Smooth orders in $\wis{alg@n}$ with center $R = \C [\C^d/G]$. Last time we have seen that Smooth orders can be described and classified locally in a combinatorial way but also that there can be no Smooth order with center $\C[\C^d/G]$. So maybe you begin to feel that I don't know what I'm talking about.

Fine, but give me one last chance to show that the overall strategy may still have some value in the desingularization project of quotient singularities. What we will see today is that there are orders $A$ over $R$ which may not be Smooth but are Smooth on a sufficiently large Zariski open subset of $\wis{rep}_{\alpha}~A$. Here 'sufficiently large' means determined by a {\em stability structure}. Whenever this is the case we can apply the results of last time to construct nice (partial) desingularizations of $\C^d/G$ and if you are in for non-commutative geometry, even a genuine non-commutative desingularization.

\section{Quotient singularities}

Last time we associated to a combinatorial triple $(Q^*,\alpha,\beta)$ a Smooth order $A^{\alpha,\beta}_{Q^*}$ with center the ring of polynomial quiver-invariants $R^{\alpha}_{Q^*}$. As we were able to classify the quiver-invariants it followed that there is no triple such that the center of $A^{\alpha,\beta}_{Q^*}$ is the coordinate ring $R = \C[\C^d/G]$ of the quotient singularity. However, there {\em are} nice orders of the form
\[
A = \frac{A^{\alpha,\beta}_{Q^*}}{I} \]
for some ideal $I$ of relations which do have center $R$ are have been used in studying quotient singularities.

\begin{mexample}[Kleinian singularities] \label{Kleinian} For a Kleinian singularity, that is, a quotient singularity  $\C^2/G$ with $G \subset SL_2(\C)$ there is an extended Dynkin diagram $D$ associated. 

Let $Q$ be the {\em double quiver} of $D$, that is to each arrow $\xymatrix{\vtx{} \ar[r]^x & \vtx{}}$ in $D$ we adjoin an arrow $\xymatrix{\vtx{} & \vtx{} \ar[l]_{x^*} }$ in $Q$ in the opposite direction and let $\alpha$ be the unique  minimal dimension vector such that $\chi_D(\alpha,\alpha) = 0$. Further, consider the {\em moment element}
\[
m = \sum_{x \in D} [x,x^*] \]
in the order $A^{\alpha}_Q$ then 
\[
A = \frac{A^{\alpha}_Q}{( m )}  \]
is an order with center $R = \C[\C^2/G]$ which is isomorphic to the skew-group algebra $\C[x,y] \# G$. Moreover, $A$ is Morita equivalent to the {\em preprojective algebra} which is the quotient of the path algebra of $Q$ by the ideal generated by the moment element
\[
\Pi_0 = \C Q/ (\sum [x,x^*] ) \]
For more details we refer to the lecture notes by W. Crawley-Boevey \cite{CrawleyLectNotes}.
\end{mexample}

\begin{mexample} \label{generalcase} Consider a quotient singularity $X = \C^d/G$ with $G \subset SL_d(\C)$ and $Q$ be the {\em McKay quiver} of $G$ acting on $V=\C^d$. 

That is, the vertices $\{ v_1,\hdots,v_k \}$ of $Q$ are in one-to-one correspondence with the irreducible representations $\{ R_1,\hdots,R_k \}$ of $G$ such that  $R_1 = \C_{triv}$ is the trivial representation. Decompose the tensorproduct in irreducibles
\[
V \otimes_{\C} R_j = R_1^{\oplus j_1} \oplus \hdots \oplus R_k^{\oplus j_k} \]
then the number of arrows in $Q$ from $v_i$ to $v_j$
\[
\#~(v_i \rTo v_j ) = j_i \]
is the multiplicity of $R_i$ in $V \otimes R_j$. Let $\alpha = (e_1,\hdots,e_k)$ be the dimension vector where $e_i = dim_{\C}~R_i$. 

The relevance of this quiver-setting is that
\[
\wis{rep}_{\alpha}~Q = Hom_G(R,R \otimes V) \]
where $R$ is the {\em regular representation}, see for example \cite{CrawNotes}. Consider $Y \subset \wis{rep}_{\alpha}~Q$ the affine subvariety of all $\alpha$-dimensional representations of $Q$ for which the corresponding $G$-equivariant map $B \in Hom_G(R,V \otimes R)$ satisfies
\[
B \wedge B = 0 \in Hom_G(R,\wedge^2 V \otimes R) \]
$Y$ is called the {\em variety of commuting matrices} and its defining relations can be expressed as linear equations between paths in $Q$ evaluated in $\wis{rep}_{\alpha}~Q$, say $(l_1,\hdots,l_z)$. Then,
\[
A = \frac{A^{\alpha}_Q}{(l_1,\hdots,l_z)} \]
is an order with center $R = \C[\C^d/G]$. In fact, $A$ is just the skew group algebra
\[
A = \C[x_1,\hdots,x_d] \# G \]
\end{mexample}

Let us give one explicit example illustrating both approaches to the Kleinian singularity $\C^2/\Z_3$.

\begin{mexample} Consider the natural action of $\Z_3$ on $\C^2$ via its embedding in $SL_2(\C)$ sending the generator to the matrix
\[
\begin{bmatrix} 
\rho & 0 \\ 0 & \rho^{-1} \end{bmatrix} 
\]
where $\rho$ is a primitive $3$-rd root of unity. $\Z_3$ has three one-dimensional simples $R_1 = \C_{triv}, R_2 = \C_{\rho}$ and $R_2 = \C_{\rho^2}$. As $V = \C^2 = R_2 \oplus R_3$ it follows that the McKay quiver setting $(Q,\alpha)$ is
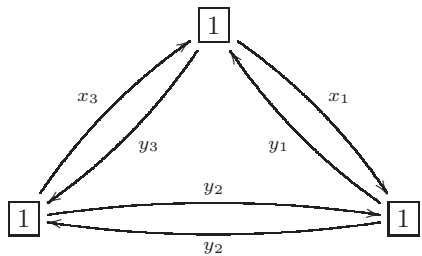
\[
\xymatrix{
& & \vtx{1} \ar@/^/[lldd]^{y_3} \ar@/^/[rrdd]^{x_1} &&  \\
& & & & \\
\vtx{1} \ar@/^/[rruu]^{x_3} \ar@/^/[rrrr]^{y_2} & & & & \vtx{1} \ar@/^/[lluu]^{y_1} \ar@/^/[llll]^{y_2} }
\]
Consider the matrices
\[
X = \begin{bmatrix} 0 & 0 & x_3 \\ x_1 & 0 & 0 \\ 0 & x_2 & 0 \end{bmatrix} \qquad \text{and} \qquad
Y = \begin{bmatrix} 0 & y_1 & 0 \\ 0 & 0 & y_2 \\ y_3 & 0 & 0 \end{bmatrix} \]
then the variety of commuting matrices is determined by the matrix-entries of $[X,Y]$ that is
\[
I = (x_3y_3-y_1x_1,x_1y_1-y_2x_2,x_2y_2-y_3x_3) \]
so the skew-group algebra is the quotient of the Smooth order $A^{\alpha}_{Q}$ (which incidentally is one of our zero-settings for dimension $4$)
\[
\C[x,y] \# \Z_3 \simeq \frac{A^{\alpha}_Q}{(x_3y_3-y_1x_1,x_1y_1-y_2x_2,x_2y_2-y_3x_3)}
\]
Taking $y_i = x_i^*$ this coincides with the description via preprojective algebras as the moment element is
\[
m = \sum_{i=1}^3 [x_i,x_i^*] = (x_3y_3-y_1x_1) e_1 + (x_1y_1-y_2x_2)e_2 + (x_2y_2-y_3x_3)e_3 \]
where the $e_i$ are the vertex-idempotents.
\end{mexample}

\begin{boxitpara}{box 0.85 setgray fill}
{\bf jotter  : }

Many interesting examples of orders are of the following form : 
\[
A = \frac{A^{\alpha}_{Q^*}}{I} \]
satisfying the following conditions :
\begin{itemize}
\item{$\alpha = (e_1,\hdots,e_k)$ is the dimension vector of a simple representation of $A$, and}
\item{the center $R = Z(A)$ is an integrally closed domain.}
\end{itemize}
These requirements (which are often hard to verify!) imply that $A$ is an order over $R$ in $\wis{alg@n}$ where $n$ is the total dimension of the simple representation, that is $| \alpha | = \sum_i e_i$.
\end{boxitpara}

Observe that such orders occur in the study of quotient singularities (see above) or as the \'etale local structure of (almost all) orders. From now on, this will be the setting we will work in.

\section{Stability structures}

For $A = A^{\alpha}_{Q^*}/I$ we define the affine variety of $\alpha$-dimensional representations
\[
\wis{rep}_{\alpha}~A = \{ V \in \wis{rep}_{\alpha}~Q^*~|~r(V)=0~\forall r \in I \}
\]
The action of $GL(\alpha) = \prod_i GL_{e_i}$ by basechange on $\wis{rep}_{\alpha}~Q^*$ induces an action (actually of $PGL(\alpha)$) on $\wis{rep}_{\alpha}~A$. Usually, $\wis{rep}_{\alpha}~A$ will have singularities but it may be smooth on the Zariski open subset 
of $\theta$-semistable representations which we will now define.

A {\em character} of $GL(\alpha)$ is determined by an integral $k$-tuple $\theta = (t_1,\hdots,t_k) \in \Z^k$
\[
\chi_{\theta}~:~GL(\alpha) \rTo \C^* \qquad (g_1,\hdots,g_k) \mapsto det(g_1)^{t_1} \hdots det(g_k)^{t_k} \]
Characters define {\em stability structures} on $A$-representations but as the acting group on $\wis{rep}_{\alpha}~A$ is really $PGL(\alpha) = GL(\alpha)/\C^*(1_{e_1},\hdots,1_{e_k})$ we only consider characters $\theta$ satisfying $\theta.\alpha = \sum_i t_ie_i = 0$.

If $V \in \wis{rep}_{\alpha}~A$ and $V' \subset V$ is an $A$-subrepresentation, that is $V' \subset V$ as representations of $Q^*$ and in addition $I(V') = 0$, we denote the dimension vector of $V'$ by $\wis{dim} V'$.

\begin{definition} For $\theta$ satisfying $\theta.\alpha = 0$, a representation $V \in \wis{rep}_{\alpha}~A$ is said to be
\begin{itemize}
\item{{\em $\theta$-semistable} if and only if for every proper $A$-subrepresentation $0 \not= V' \subset V$ we have $\theta.\wis{dim} V' \geq 0$.}
\item{{\em $\theta$-stable} if and only if for every proper $A$-subrepresentation $0 \not= V' \subset V$ we have $\theta.\wis{dim} V' > 0$.}
\end{itemize}
\end{definition}

For any setting $\theta.\alpha = 0$ we have the following inclusions of Zariski open $GL(\alpha)$-stable subsets of $\wis{rep}_{\alpha}~A$
\[
\wis{rep}^{simple}_{\alpha}~A \subset \wis{rep}^{\theta-stable}_{\alpha}~A \subset \wis{rep}^{\theta-semist}_{\alpha}~A \subset \wis{rep}_{\alpha}~A \]
but one should note that some of these open subsets may actually be empty!

Recall that a point of the algebraic quotient variety $\wis{iss}_{\alpha}~A = \wis{rep}_{\alpha} // GL(\alpha)$ represents the orbit of an $\alpha$-dimensional semi-simple representation $V$ and such representations can be separated by the values $f(V)$ where $f$ is a polynomial invariant on $\wis{rep}_{\alpha}~A$. This follows  because the coordinate ring of the quotient variety
\[
\C[\wis{iss}_{\alpha}~A] = \C[\wis{rep}_{\alpha}~A]^{GL(\alpha)} \]
and points correspond to maximal ideals of this ring. Recall from \cite{LBProcesi} that the invariant ring is generated by taking traces along oriented cycles in the marked quiver-setting $(Q^*,\alpha)$.

\begin{boxitpara}{box 0.85 setgray fill}
{\bf jotter  : }

For $\theta$-stable and $\theta$-semistable representations there are similar results and morally one should view $\theta$-stable representations as corresponding to simple representations whereas $\theta$-semistables are arbitrary representations. 

For this reason we will only be able to classify direct sums of $\theta$-stable representations by certain algebraic varieties which are called the {\em moduli spaces} of semistables representations.
\end{boxitpara}

The notion corresponding to a polynomial invariant in this more general setting is that of a {\em polynomial semi-invariant}. A polynomial function $f \in \C[\wis{rep}_{\alpha}~A]$ is said to be a $\theta$-semi-invariant of {\em weight} $l$ if for all $g \in GL(\alpha)$ we have
\[
g.f = \chi_{\theta}(g)^l f \]
where $\chi_{\theta}$ is the character of $GL(\alpha)$ corresponding to $\theta$. A representation $V \in \wis{rep}_{\alpha}~A$ is $\theta$-semistable if and only if there is a $\theta$-semi-invariant $f$ of some weight $l$ such that $f(V) \not= 0$.

Clearly, $\theta$-semi-invariants of weight zero are just polynomial invariants and the multiplication of $\theta$-semi-invariants of weight $l$ resp. $l'$ has weight $l+l'$. Hence, the ring of all $\theta$-semi-invariants
\[
\C[\wis{rep}_{\alpha}~A]^{GL(\alpha),\theta} = \oplus_{l=0}^{\infty} \{ f \in \C[\wis{rep}_{\alpha}~A]~| \forall g \in GL(\alpha)~:~g.f = \chi_{\theta}(g)^l f~\} \]
is a graded algebra with part of degree zero $\C[\wis{iss}_{\alpha}~A]$. But then we have a {\em projective morphism}
\[
\wis{proj}~\C[\wis{rep}_{\alpha}~A]^{GL(\alpha),\theta} \rOnto^{\pi} \wis{iss}_{\alpha}~A \]
such that all fibers of $\pi$ are projective varieties. The main properties of $\pi$ can be deduced from  \cite{King}

\begin{theorem} Points in $\wis{proj}~\C[\wis{rep}_{\alpha}~A]^{GL(\alpha),\theta}$ are in one-to-one correspondence with isomorphism classes of direct sums of $\theta$-stable representations of total dimension $\alpha$.

If $\alpha$ is such that there are $\alpha$-dimensional simple $A$-representations, then $\pi$ is a birational map. 
\end{theorem}

\begin{definition}
We call $\wis{proj}~\C[\wis{rep}_{\alpha}~A]^{GL(\alpha),\theta}$ the {\em moduli space} of $\theta$-semistable representations of $A$ and denote it with $\wis{moduli}^{\theta}_{\alpha}~A$.
\end{definition}

\begin{mexample} In the case of Kleinian singularities, see example~\ref{Kleinian}, if we take $\theta$ to be a generic character such that $\theta.\alpha = 0$, then the projective map
\[
\wis{moduli}^{\theta}_{\alpha}~A \rOnto X = \C^2/G \]
is a minimal resolution of singularities. Note that the map is birational as $\alpha$ is the dimension vector of a simple representation of $A = \Pi_0$, see \cite{CrawleyLectNotes}.
\end{mexample}

\begin{mexample} For general quotient singularities, see example~\ref{generalcase}, assume that the first vertex in the McKay quiver corresponds to the trivial representation. Take a character $\theta \in \Z^k$ such that $t_1 < 0$ and all $t_i > 0$ for $i \geq 2$, for example take
\[
\theta = ( - \sum_{i=2}^k \wis{dim} R_i , 1, \hdots, 1 ) \]
Then, the corresponding moduli space is isomorphic to
\[
\wis{moduli}^{\theta}_{\alpha}~A \simeq G-\wis{Hilb}~\C^d \]
the {\em $G$-equivariant Hilbert scheme} which classifies all $\# G$-codimensional ideals $I \triangleleft \C[x_1,\hdots,x_d]$ where
\[
\frac{\C[x_1,\hdots,x_d]}{I} \simeq \C G \]
as $G$-modules, hence in particular $I$ must be stable under the action of $G$. It is well known that the natural map
\[
G-\wis{Hilb}~\C^d \rOnto X = \C^d/G \]
is a minimal resolution if $d=2$ and if $d=3$ it is often a crepant resolution, for example whenever $G$ is Abelian. In non-Abelian cases it may have remaining singularities though which often are of conifold type. See \cite{CrawNotes} for more details.
\end{mexample}

\begin{boxitpara}{box 0.85 setgray fill}
{\bf jotter  : }

My motivation for this series of talks was to look for a non-commutative explanation for the omnipresence of conifold singularities in partial resolutions of three dimensional quotient singularities as well as to have a conjectural list of possible remaining singularities for higher dimensional quotient singularities. 
\end{boxitpara}

\begin{mexample}  In the $\C^2/\Z_3$-example one can take $\theta=(-2,1,1)$. The following representations
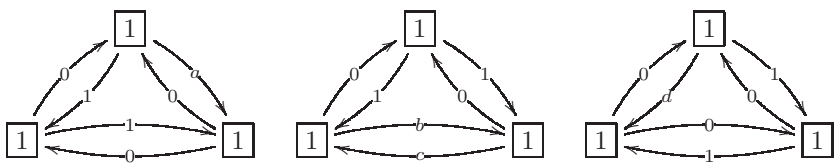
\[
\xymatrix@=.3cm{
& & \vtx{1} \ar@/^/[lldd]|{1} \ar@/^/[rrdd]|{a} &&  \\
& & & & \\
\vtx{1} \ar@/^/[rruu]|{0} \ar@/^/[rrrr]|{1} & & & & \vtx{1} \ar@/^/[lluu]|{0} \ar@/^/[llll]|{0} } \quad
\xymatrix@=.3cm{
& & \vtx{1} \ar@/^/[lldd]|{1} \ar@/^/[rrdd]|{1} &&  \\
& & & & \\
\vtx{1} \ar@/^/[rruu]|{0} \ar@/^/[rrrr]|{b} & & & & \vtx{1} \ar@/^/[lluu]|{0} \ar@/^/[llll]|{c} } \quad
\xymatrix@=.3cm{
& & \vtx{1} \ar@/^/[lldd]|{d} \ar@/^/[rrdd]|{1} &&  \\
& & & & \\
\vtx{1} \ar@/^/[rruu]|{0} \ar@/^/[rrrr]|{0} & & & & \vtx{1} \ar@/^/[lluu]|{0} \ar@/^/[llll]|{1} } \]
are all nilpotent and are $\theta$-stable. In fact if $bc=0$ they are representants of the exceptional fiber of the desingularization
\[
\wis{moduli}^{\theta}_{\alpha}~A \rOnto \wis{iss}_{\alpha}~A = \C^2/\Z_3 \]
\end{mexample}

\section{Partial resolutions}

It is about time we state the main result of these notes which was proved in \cite{LBSymens}.

\begin{theorem} Let $A = A^{\alpha}_{Q^*}/(R)$ be an $R$-order in $\wis{alg@n}$. Assume that there exists a stability structure $\theta \in \Z^k$ such that the Zariski open subset $\wis{rep}^{\theta-semist}_{\alpha}~A$ of all $\theta$-semistable $\alpha$-dimensional representations of $A$ is a smooth variety. 

Then there exists a sheaf $\Ascr$ of Smooth orders over $\wis{moduli}^{\theta}_{\alpha}~A$ such that the diagram below is commutative
\[
\begin{diagram}
\wis{spec}~\Ascr & & \\
\dTo^{c} & \rdTo^{\phi} \\
\wis{moduli}^{\theta}_{\alpha}~A & \rOnto^{\pi} & X = \wis{spec}~R
\end{diagram}
\]
Here, $\wis{spec}~\Ascr$ is a non-commutative variety obtained by gluing affine non-commutative varieties $\wis{spec}~A_i$ together and $c$ is the map which intersects locally a maximal ideal with the center. As $\Ascr$ is a sheaf of Smooth orders, $\phi$ can be viewed as a non-commutative desingularization of $X$.

If you are only interested in commutative desingularizations,  $\pi$ is a partial resolution of $X$ and we have full control over the remaining singularities in $\wis{moduli}^{\theta}_{\alpha}~A$, that is, all remaining singularities are of the form classified in the previous lecture.

Moreover, if $\theta$ is such that all $\theta$-semistable $A$-representations are actually $\theta$-stable, then $\Ascr$ is a sheaf of Azumaya algebras over $\wis{moduli}^{\theta}_{\alpha}~A$ and in this case $\pi$ is a commutative desingularization of $X$. If, in addition, also $gcd(\alpha) = 1$, then $\Ascr \simeq End~\Pscr$ for some vectorbundle of rank $n$ over $\wis{moduli}^{\theta}_{\alpha}~A$.
\end{theorem}

\begin{boxitpara}{box 0.85 setgray fill}
{\bf jotter  : }

It should be stressed that the condition that $\wis{rep}^{\theta-semist}_{\alpha}~A$ is a smooth variety is {\em very strong} and is usually {\em hard} to verify in concrete situations.
\end{boxitpara}

\begin{mexample} In the case of Kleinian singularities, see example~\ref{Kleinian}, there exists a suitable stability structure $\theta$ such that $\wis{rep}^{\theta-semist}_{\alpha}~\Pi_0$ is smooth.
For consider the {\em moment map}
\[
\wis{rep}_{\alpha}~Q \rTo^{\mu} \wis{lie}~GL(\alpha) = M_{\alpha}(\C) = M_{e_1}(\C) \oplus \hdots \oplus M_{e_k}(\C)
\]
defined by sending $V = (V_a,V_{a^*})$ to
\[
 (\sum_{\xymatrix{\vtx{} \ar[r]^a&\vtx{1}}} V_aV_{a^*} - \sum_{\xymatrix{\vtx{1} \ar[r]^a & \vtx{}}} V_{a^*}V_a, \hdots, \sum_{\xymatrix{\vtx{} \ar[r]^a&\vtx{k}} }V_aV_{a^*} - \sum_{\xymatrix{\vtx{k} \ar[r]^a & \vtx{}} }V_{a^*}V_a) \]
The differential $d \mu$ can be verified to be surjective in any representation $V \in \wis{rep}_{\alpha}~Q$ which has stabilizer subgroup $\C^*(1_{e_1},\hdots,1_{e_k})$ (a so called {\em Schur representation}) see for example \cite[lemma 6.5]{CrawleyMoment}. 

Further, any $\theta$-stable representation is Schurian. Moreover, for a generic stability structure $\theta \in \Z^k$ we have that every $\theta$-semistable $\alpha$-dimensional representation is $\theta$-stable as the $gcd(\alpha) = 1$.
Combining these facts it follows that $\mu^{-1}(0) = \wis{rep}_{\alpha}~\Pi_0$ is smooth in all $\theta$-stable representations.
\end{mexample}

\begin{mexample} Another case where smoothness of $\wis{rep}^{\theta-semist}_{\alpha}~A$ is evident is when $A = A^{\alpha}_{Q^*}$ is a Smooth order as then $\wis{rep}_{\alpha}~A$ itself is smooth. This observation can be used to resolve the remaining singularities in the partial resolution.

If $gcd(\alpha) = 1$ then for a sufficiently general $\theta$ all $\theta$-semistable representations are actually $\theta$-stable whence the quotient map
\[
\wis{rep}^{\theta-semist}_{\alpha}~A \rOnto \wis{moduli}^{\theta}_{\alpha}~A \]
is a principal $PGL(\alpha)$-fibration and as the total space is smooth, so is $\wis{moduli}^{\theta}_{\alpha}~A$. Therefore, the projective map
\[
\wis{moduli}^{\theta}_{\alpha}~A \rOnto^{\pi} \wis{iss}_{\alpha}~A \]
is a resolution of singularities in this case.

However, if $l = gcd(\alpha)$, then $\wis{moduli}_{\alpha}^{\theta}~A$ will usually contain singularities which are as bad as the quotient variety singularity of tuples of $l \times l$ matrices under simultaneous conjugation.
\end{mexample}

Fortunately, the proof of the theorem will follow from the hard work we did in last lecture provided we can solve two problems. 

A minor problem is that  we classified central singularities of Smooth orders in $\wis{alg@n}$ but here we are working with $\alpha$-dimensional representations and with the action of $GL(\alpha)$ rather than $GL_n$. This problem we will address immediately.

A more serious problem is that $\wis{rep}^{\theta-semist}_{\alpha}~A$ is not an affine variety in general so we will have to cover it with affine varieties $X_i$ and consider associated orders $A_i$. But then we have to clarify why $\theta$-semistable representations of $A$ correspond to {\em all} representations of the $A_i$. This may not be clear at first sight.

\section{Going from $\wis{alg@n}$ to $\wis{alg@}\alpha$}

If $Q^*$ is a marked quiver on $k$ vertices, then the subalgebra generated by the vertex-idempotents $\C^k$ is a subalgebra of $A = A^{\alpha}_{Q^*}/(R)$ hence we have a morphism
\[
\wis{rep}_n~A \rTo \wis{rep}_n~\C^k = \bigsqcup_{| \alpha | = n} GL_n/GL(\alpha) \]
where the last decomposition follows from the fact that $\C^k$ is semi-simple whence every $n$-dimensional representation is fully determined by the multiplicities of the simple $1$-dimensional components. 

Further, we should consider $\wis{trep}_n~A$ the subvariety of trace preserving $A$-representations but a trace map on $A$ fixes the trace on $\C^k$ and hence determines the component $GL_n/GL(\alpha)$. That is, we have that
\[
\wis{trep}_n~A = GL_n \times^{GL(\alpha)} \wis{rep}_{\alpha}~A \]
the variety is a {\em principal fiber bundle}. 

That is, if $V$ is any $n$-dimensional trace preserving $A$-representation $A \rTo^{\phi} M_n(\C)$ then the images $\phi(v_i)$ of the vertex-idempotents are a full set of orthogonal idempotents so they can be conjugated to a set of matrices 
\[
\phi'(v_i) = \begin{bmatrix} \ddots & & & & \\
& 1 & & & \\
& & \ddots & & \\
& & & 1 & \\
& & & & \ddots \end{bmatrix}
\]
with only $1$'s from place $\sum_{j=1}^{i-1} e_j + 1$ to place $\sum_{j=1}^i e_j$. But using these idempotents we see that the representation $\phi'~:~A \rTo M_n(\C)$ has block-matrices coming from a representation in $\wis{rep}_{\alpha}~A$.

As is the case for any principal fiber bundle, this gives a natural one-to-one correspondence between
\begin{itemize}
\item{$GL_n$-orbits in $\wis{trep}_n~A$, and}
\item{$GL(\alpha)$-orbits in $\wis{rep}_{\alpha}~A$.}
\end{itemize}
Moreover the corresponding quotient varieties $\wis{tiss}_n~A = \wis{trep}_n~A // GL_n$ and $\wis{iss}_{\alpha}~A = \wis{rep}_{\alpha}~A // GL(\alpha)$ are isomorphic so we can apply all our $(P)GL_n$-results to this setting.

\begin{boxitpara}{box 0.85 setgray fill}
{\bf jotter  : }

Alternatively, we can define $\wis{alg@}\alpha$ to be the subcategory of $\wis{alg@n}$ with objects the algebras $A \in \wis{alg@n}$ which are $\C^k$-algebras via the embedding given by the matrices $\phi'(v_i)$ above and with morphism the $\C^k$-algebra morphisms in $\wis{alg@n}$. 

It is then clear that a Smooth order in $\wis{alg@}\alpha$ (that is, having the lifting property with respect to nilpotent ideals in $\wis{alg@}\alpha$) is a Smooth order in $\wis{alg@n}$ which is an object in  $\wis{alg@}\alpha$.
\end{boxitpara}

\section{The affine opens  $X_D$}

To solve the second problem, we claim that we can cover the moduli space 
\[
\wis{moduli}^{\theta}_{\alpha}~A = \bigcup_D X_D \]
where $X_D$ is an affine open subset such that under the canonical quotient map
\[
\wis{rep}^{\theta-semist}_{\alpha}~A \rOnto^{\pi} \wis{moduli}^{\theta}_{\alpha}~A \]
we have that
\[
\pi^{-1}(X_D) = \wis{rep}_{\alpha}~A_D \]
for some $\C[X_D]$-order $A_D$ in $\wis{alg@n}$. 

If in addition $\wis{rep}_{\alpha}^{\theta-semist}~A$ is a smooth variety, each of the $\wis{rep}_{\alpha}~A_D$ are smooth affine $GL(\alpha)$-varieties whence the orders $A_D$ are all Smooth and the result will follow from the results of last lecture.

Because $\wis{moduli}^{\theta}_{\alpha}~A = \wis{proj}~\C[\wis{rep}_{\alpha}~A]^{GL(\alpha),\theta}$ we need control on the generators of all $\theta$-semi-invariants. Such a generating set was found by Aidan Schofield and Michel Van den Bergh in \cite{SchofVdB} : {\em determinantal semi-invariants}. In order to define them we have to introduce some notation first.

Reorder the vertices in $Q^*$ such that the entries of $\theta$ are separated in three strings
\[
\theta = (\underbrace{t_1,\hdots,t_i}_{> 0},\underbrace{t_{i+1},\hdots,t_j}_{=0},\underbrace{t_{j+1},\hdots,t_k}_{< 0}) \]
and let $\theta$ be such that $\theta.\alpha = 0$. Fix a weight $l \in \N_+$ and take arbitrary natural numbers $\{ l_{i+1},\hdots,l_j \}$. 

Consider a rectangular matrix $L$ with
\begin{itemize}
\item{$lt_1+\hdots+lt_i+l_{i+1} + \hdots + l_j$ rows and}
\item{$l_{i+1} + \hdots + l_j - l t_{j+1} - \hdots - l t_k$ columns}
\end{itemize}
\[
L = \quad \begin{array}{cc||c|c|c|c|c|c}
& & \overbrace{}^{l_{i+1}} & \hdots & \overbrace{}^{l_j} & \overbrace{}^{-lt_{j+1}} & \hdots & \overbrace{}^{-lt_k} \\
\hline \hline 
lt_1&  \{ & L_{1,i+1} & & L_{1,j} & L_{1,j+1} & & L_{1,k} \\
\hline
& \vdots & & & & & \\
\hline 
lt_i & \{ & L_{i,i+1} & & L_{i,j} & L_{i,j+1} & & L_{i,k} \\
\hline
l_{i+1} & \{ & L_{i+1,i+1} & & L_{i+1,j} & L_{i+1,j+1} & & L_{i+1,k} \\
\hline
& \vdots & & & & & \\
\hline
l_j & \{ & L_{j,i+1} & & L_{j,j}  & L_{j,j+1} & & L_{j,k} 
\end{array}
\]

in which each entry of $L_{r,c}$ is a linear combination of oriented paths in the marked quiver $Q^*$ with starting vertex $v_c$ and ending vertex $v_r$.

The relevance of this is that we can evaluate $L$ at any representation $V \in \wis{rep}_{\alpha}~A$ and obtain a {\em square matrix} $L(V)$ as $\theta.\alpha = 0$. More precisely, if $V_i$ is the vertex-space of $V$ at vertex $v_i$ (that is, $V_i$ has dimension $e_i$), then evaluating $L$ at $V$ gives a linear map
\[
\begin{diagram}
V_{i+1}^{\oplus l_{i+1}} \oplus \hdots \oplus V_j^{\oplus l_j} \oplus V_{j+1}^{\oplus -lt_{j+1}} \oplus \hdots \oplus V_k^{\oplus -lt_k} \\
\dTo^{L(V)} \\
V_1^{\oplus lt_1} \oplus \hdots \oplus V_i^{\oplus lt_i} \oplus V_{i+1}^{\oplus l_{i+1}} \oplus \hdots \oplus V_j^{\oplus l_j}
\end{diagram}
\]
and $L(V)$ is a square $N \times N$ matrix where
\[
l_{i+1} + \hdots + l_j - lt_{j+1} - \hdots - lt_k = N = lt_1 + \hdots + lt_i + l_{i+1} + \hdots + l_j \]
So we can consider $D(V) = \wis{det} L(V)$ and verify that $D$ is a $GL(\alpha)$-semi-invariant polynomial on $\wis{rep}_{\alpha}~A$ of weight $\chi_{\theta}^l$. The result of \cite{SchofVdB} asserts that these {\em determinantal semi-invariants} are algebra generators of the graded algebra
\[
\C[\wis{rep}_{\alpha}~A]^{GL(\alpha),\theta} \]
Observe that this result is to semi-invariants what the result of \cite{LBProcesi} is to invariants. In fact, one can deduce the latter from the first. 

We have seen that a representation $V \in \wis{rep}_{\alpha}~A$ is $\theta$-semistable if and only if some semi-invariant of weight $\chi_{\theta}^l$ for some $l$ is non-zero on it. This proves

\begin{theorem} The Zariski open subset of $\theta$-semistable $\alpha$-dimensional $A$-representations can be covered by affine $GL(\alpha)$-stable open subsets
\[
\wis{rep}^{\theta-semist}_{\alpha}~A = \bigsqcup_D \{ V~|~D(V) = \wis{det} L(V) \not= 0 \} \]
and hence the moduli space can also be covered by affine open subsets
\[
\wis{moduli}^{\theta}_{\alpha}~A = \bigcup_D~X_D \]
where
$
X_D = \{ [V] \in \wis{moduli}^{\theta}_{\alpha}~A~|~D(V)=\wis{det} L(V) \not= 0 \}
$.
\end{theorem}

\begin{mexample} \label{voorbeeld} In the $\C^2/\Z_3$ example, the $\theta$-semistable representations
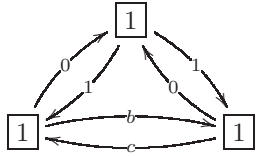
\[
\xymatrix@=.3cm{
& & \vtx{1} \ar@/^/[lldd]|{1} \ar@/^/[rrdd]|{1} &&  \\
& & & & \\
\vtx{1} \ar@/^/[rruu]|{0} \ar@/^/[rrrr]|{b} & & & & \vtx{1} \ar@/^/[lluu]|{0} \ar@/^/[llll]|{c} }
\]
with $\theta=(-2,1,1)$ all lie in the affine open subset $X_D$ where $L$ is a matrix of the form
\[
L = \begin{bmatrix} x_1 & 0 \\ \ast & y_3 \end{bmatrix} \]
where $\ast$ is any path in $Q$ starting in $x_1$ and ending in $x_3$.
\end{mexample}

\section{The $\C[X_D]$- orders $A_D$}

Analogous to the rectangular matrix $L$ we define a rectangular matrix $N$ with
\begin{itemize}
\item{$lt_1+\hdots+lt_i+l_{i+1} + \hdots + l_j$ columns and}
\item{$l_{i+1} + \hdots + l_j - l t_{j+1} - \hdots - l t_k$ rows}
\end{itemize}
\[
N = \quad \begin{array}{cc||c|c|c|c|c|c}
& & \overbrace{}^{l t_1} & \hdots & \overbrace{}^{l t_i} & \overbrace{}^{l_{i+1}} & \hdots & \overbrace{}^{l_j} \\
\hline \hline
l_{i+1} &  \{ & N_{i+1,1} & & N_{i+1,i} & N_{i+1,i+1} & & N_{i+1,j} \\
\hline
& \vdots & & & & & \\
\hline
l_j & \{ & N_{j,1} & & N_{j,i} & N_{j,i+1} & & N_{j,j} \\
\hline
-lt_{j+1}  & \{ & N_{j+1,1} & & N_{j+1,i} & N_{j+1,i+1} & & N_{j+1,j} \\
\hline
& \vdots & & & & & \\
\hline
-l t_k & \{ & N_{k,1} & & N_{k,i}  & N_{k,i+1} & & N_{k,j} 
\end{array}
\]
filled with new variables and define an {\em extended marked quiver} $Q^*_D$ where we adjoin for each entry in $N_{r,c}$ an additional arrow from $v_c$ to $v_r$ and denote it with the corresponding variable from $N$.

Let $I_1$ (resp. $I_2$ be the set of relations in $\C Q^*_D$ determined from the matrix-equations
{\tiny
\[
N.L = \begin{bmatrix} \boxed{(v_{i+1})_{l_{i+1} } } & & & & & 0 \\
& \ddots & & & & \\
& & \boxed{(v_j)_{l_j}} & & & \\
& & & \boxed{(v_{j+1})_{-lt_{j+1}} }& &  \\
& & & & \ddots & \\
0 & & & & & \boxed{(v_k)_{-lt_k}}
\end{bmatrix}
\]}
respectively
{\tiny
\[
L.N = \begin{bmatrix}
\boxed{(v_1)_{lt_1}} & & & & & 0 \\
& \ddots & & & & \\
& & \boxed{(v_i)_{lt_i} }& & & \\
& & & \boxed{(v_{i+1})_{l_{i+1}} }& & \\
& & & & \ddots & \\
0 & & & & & \boxed{(v_j)_{l_j}}
\end{bmatrix}
\]}
where $(v_i)_{n_j}$ is the square $n_j \times n_j$ matrix with $v_i$ on the diagonal and zeroes elsewhere.

Define a new non-commutative order
\[
A_D = \frac{A^{\alpha}_{Q^*_D}}{(I,I_1,I_2)}
\]
then $A_D$ is a $\C[X_D]$-order in $\wis{alg@n}$. 

\begin{mexample}
In the setting of example~\ref{voorbeeld} with $\ast = y_3$, the extended quiver-setting $(Q_D,\alpha)$ is
\[
\xymatrix@=.3cm{
& & \vtx{1} \ar@/^/[lldd]|{y_3} \ar@/^/[rrdd]|{x_1} &&  \\
& & & & \\
\vtx{1} \ar@/^/[rruu]|{x_3} \ar@/^/[rrrr]|{y_2} \ar@/^4ex/[rruu]|{n_3} \ar@/^6ex/[rruu]|{n_4} & & & & \vtx{1} \ar@/^/[lluu]|{y_1} \ar@/^/[llll]|{x_2} \ar@/_4ex/[lluu]|{n_1} \ar@/_6ex/[lluu]|{n_2} }
\]
Hence, with
\[
L = \begin{bmatrix} x_1 & 0 \\ y_3 & y_3 \end{bmatrix} \qquad N = \begin{bmatrix} n_1 & n_3 \\ n_2 & n_4 \end{bmatrix} \]
the defining equations of the order $A_D$ become
\[
\begin{cases} I = (x_3y_3-y_1x_1,x_1y_1-y_2x_2,x_2y_2-y_3x_3) \\
I_1 = (n_1x_1+n_3y_3-v_1,n_3y_3,n_2x_1+n_4y_3,n_4y_3-v_1) \\
I_2 = (x_1n_1-v_2,x_1n_3,y_3n_1+y_2n_2,y_3n_3+y_3n_4-v_3)
\end{cases}
\]
\end{mexample}

\begin{boxitpara}{box 0.85 setgray fill}
{\bf jotter  : }

This construction may seem a bit mysterious at first  but what we are really doing is to construct the {\em universal localization} as in for example \cite{Schofield} associated to the map between projective $A$-modules determined by $L$, but this time not in the category $\wis{alg}$ of all algebras but in $\wis{alg@}\alpha$.
\end{boxitpara}

That is, take $P_i = v_i A$ be the projective right ideal associated to vertex $v_i$, then $L$ determines an $A$-module morphism
\[
P = P_{i+1}^{\oplus l_{i+1}} \oplus \hdots \oplus P_k^{\oplus -lt_k} \rTo^{L} P_1^{\oplus lt_1} \oplus \hdots \oplus P_j^{\oplus l_j} = Q \]
The algebra map $A \rTo^{\phi} A_D$ is universal in $\wis{alg@}\alpha$ with respect to $L \otimes \phi$ being invertible, that is, if $A \rTo^{\psi} S$ is a morphism in $\wis{alg@}\alpha$ such that $L \otimes \psi$ is an isomorphism of right $S$-modules, then there is a unique map in $\wis{alg@}\alpha$ $A_D \rTo^u S$ such that $\psi = u \circ \phi$.

The proof of the main result follows from the following result :
\[
\begin{diagram}
\wis{rep}^{\theta-semist}~A & \lInto & \pi^{-1}(X_D) \simeq \wis{rep}_{\alpha}~A_D \\
\dOnto^{\pi} & & \dOnto \\
\wis{moduli}^{\theta}_{\alpha}~A & \lInto & X_D
\end{diagram}
\]

\begin{theorem} The following statements are equivalent
\begin{enumerate}
\item{$V \in \wis{rep}_{\alpha}^{\theta-semist}~A$ lies in $\pi^{-1}(X_D)$, and}
\item{There is a unique extension $\tilde{V}$ of $V$ such that $\tilde{V} \in \wis{rep}_{\alpha}~A_D$.}
\end{enumerate}
\end{theorem}

\begin{proof}
$1 \Rightarrow 2$ : Because $L(V)$ is invertible we can take $N(V)$ to be its inverse and decompose it into blocks corresponding to the new arrows in $Q^*_D$. This then defines the unique extension $\tilde{V} \in \wis{rep}_{\alpha}~Q^*_D$ of $V$. As $\tilde{V}$ satisfies $R$ (because $V$ does) and $R_1$ and $R_2$ (because $N(V) = L(V)^{-1}$) we have that $\tilde{V} \in \wis{rep}_{\alpha}~A_D$.

$2 \Rightarrow 1$ :  Restrict $\tilde{V}$ to the arrows of $Q$ to get a $V \in \wis{rep}_{\alpha}~Q$. As $\tilde{V}$ (and hence $V$) satisfies $R$, $V \in \wis{rep}_{\alpha}~A$. Moreover, $V$ is such that $L(V)$ is invertible (this follows because $\tilde{V}$ satisfies $R_1$ and $R_2$). Hence, $D(V) \not= 0$ and because $D$ is a $\theta$-semiinvariant it follows that $V$ is an $\alpha$-dimensional $\theta$-semistable representation of $A$. An alternative method to see this is as follows. Assume that $V$ is {\em not} $\theta$-semistable and let $V' \subset V$ be a subrepresentation such that $\theta.\wis{dim} V' < 0$. Consider the restriction of the linear map $L(V)$ to the subrepresentation $V'$ and look at the commuting diagram
\[
\begin{diagram}
V_{i+1}^{'\oplus l_{i+1}} \oplus \hdots \oplus V_k^{'\oplus -lt_k} & \rTo^{L(V)|V'} & V_1^{'\oplus lt_1} \oplus \hdots \oplus V_j^{'\oplus l_j} \\
\dInto & & \dInto \\
V_{i+1}^{\oplus l_{i+1}} \oplus \hdots \oplus V_k^{\oplus -lt_k} & \rTo^{L(V)} & V_1^{\oplus lt_1} \oplus \hdots \oplus V_j^{\oplus l_j}
\end{diagram}
\]
As $\theta. \wis{dim} V' < 0$ the top-map must have a kernel which is clearly absurd as we know that $L(V)$ is invertible.
\end{proof}

\begin{mexample} In the setting of example~\ref{voorbeeld} with $\ast = y_3$ we have that the uniquely determined extension of the $A$-representation
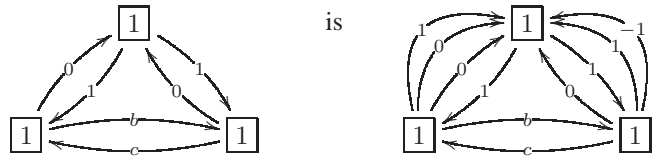
\[
\xymatrix@=.3cm{
& & \vtx{1} \ar@/^/[lldd]|{1} \ar@/^/[rrdd]|{1} &&  \\
& & & & \\
\vtx{1} \ar@/^/[rruu]|{0} \ar@/^/[rrrr]|{b} & & & & \vtx{1} \ar@/^/[lluu]|{0} \ar@/^/[llll]|{c} }~\qquad \text{is} \qquad
\xymatrix@=.3cm{
& & \vtx{1} \ar@/^/[lldd]|{1} \ar@/^/[rrdd]|{1} &&  \\
& & & & \\
\vtx{1} \ar@/^/[rruu]|{0} \ar@/^/[rrrr]|{b} \ar@/^4ex/[rruu]|{0} \ar@/^6ex/[rruu]|{1} & & & & \vtx{1} \ar@/^/[lluu]|{0} \ar@/^/[llll]|{c} \ar@/_4ex/[lluu]|{1} \ar@/_6ex/[lluu]|{-1} }
\]
Observe that this extension is a simple $A_D$-representation for every $b,c \in \C$.
\end{mexample}

\section{Non-commutative desingularizations}

There is just one more thing to clarify : how are the different $A_D$'s glued together to form a sheaf $\Ascr$ of non-commutative algebras over $\wis{moduli}^{\theta}_{\alpha}~A$ and how can we construct the non-commutative variety $\wis{spec}~\Ascr$? The solution to both problems follows from the universal property of $A_D$.

Let $A_{D_1}$ (resp. $A_{D_2}$) be the algebra constructed from a rectangular matrix $L_1$ (resp. $L_2$), then we can construct the direct sum map $L = L_1 \oplus L_2$ for which the corresponding semi-invariant $D=D_1D_2$. As $A \rTo A_D$ makes the projective module morphisms associated to $L_1$ and $L_2$ an isomorphism we have uniquely determined maps in $\wis{alg@}\alpha$
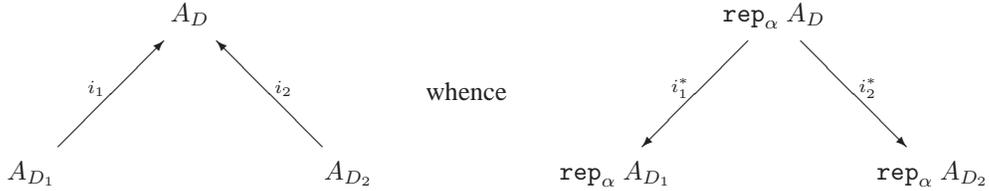
\[
\begin{diagram}
& & A_D & & \\
& \ruTo^{i_1} & & \luTo^{i_2} & \\
A_{D_1} & & & & A_{D_2}
\end{diagram}
\qquad \text{whence}
\qquad
\begin{diagram}
& & \wis{rep}_{\alpha}~A_D & & \\
& \ldTo^{i_1^*} & & \rdTo^{i^*_2} & \\
\wis{rep}_{\alpha}~A_{D_1} & & & & \wis{rep}_{\alpha}~A_{D_2} 
\end{diagram}
\]
As $\wis{rep}_{\alpha}~A_D = \pi^{-1}(X_D)$ (and similarly for $D_i$) we have that $i_j^*$ are embeddings as are the $i_j$. This way we can glue the sections $\Gamma(X_{D_1},\Ascr) = A_{D_1}$ with $\Gamma(X_{D_2},\Ascr) = A_{D_2}$ over their intersection $X_D = X_{D_1} \cap X_{D_2}$ via the inclusions $i_j$. Hence we get a coherent sheaf of non-commutative algebras $\Ascr$ over $\wis{moduli}^{\theta}_{\alpha}~A$.

Observe that many of the orders $A_D$ are isomorphic. In example~\ref{voorbeeld} all matrices $L$ with fixed diagonal entries $x_1$ and $y_3$ but with varying $\ast$-entry have isomorphic orders $A_D$ (use the universal property). 

In a similar way we would like to glue $\wis{max}~A_{D_1}$ with $\wis{max}~A_{D_2}$ over $\wis{max}~A_D$ using the algebra maps $i_j$ to form a non-commutative variety $\wis{spec}~\Ascr$. However, the construction of $\wis{max}~A$ and the non-commutative structure sheaf is {\em not} functorial in general.

\begin{mexample} Consider the inclusion map map in $\wis{alg@2}$
\[
A = \begin{bmatrix} R & R \\ I  & R \end{bmatrix} \rInto \begin{bmatrix}
R & R \\ R & R \end{bmatrix} = A' \]
then all twosided maximal ideals of $A'$ are of the form $M_2(\mathfrak{m})$ where $\mathfrak{m}$ is a maximal ideal of $R$. If $I \subset \mathfrak{m}$ then the intersection
\[
\begin{bmatrix} \mathfrak{m} & \mathfrak{m} \\ \mathfrak{m} & \mathfrak{m} \end{bmatrix} \cap \begin{bmatrix} R & R \\ I & R \end{bmatrix} = \begin{bmatrix}
\mathfrak{m} & \mathfrak{m} \\ I & \mathfrak{m} \end{bmatrix} \]
which is {\em not} a maximal ideal of $A$ as 
\[
\begin{bmatrix} \mathfrak{m} & R \\ I & R \end{bmatrix} \begin{bmatrix} R & R \\ I & \mathfrak{m} \end{bmatrix} = \begin{bmatrix} \mathfrak{m} & \mathfrak{m} \\ I & \mathfrak{m} \end{bmatrix} \]
and so there is no natural map $\wis{max} A' \rTo \wis{max}~A$, let alone a continuous one.
\end{mexample}

\begin{boxitpara}{box 0.85 setgray fill}
{\bf jotter  : }

Associating to a non-commutative algebra $A$ its prime ideal spectrum $\wis{spec}~A$ is only functorial for {\em extensions} $A \rTo^f B$, that is, satisfying
\[
B = f(A) Z_B(A) \qquad \text{with} \qquad Z_B(A) = \{ b \in B~|~bf(a) = f(a)b~\forall a \in A \} \]
In \cite{Procesibook} it was proved that if $A \rTo^f B$ is an extension then the map
\[
\wis{spec}~B \rTo \wis{spec}~R \qquad \qquad P \rTo f^{-1}(P) \]
is well-defined and continuous for the Zariski topology.
\end{boxitpara}

Fortunately, in the case of interest to us, that is for the maps $i_j~:~A_{D_j} \rTo A_D$ this presents no problem as they are even {\em central extensions}, that is
\[
A_D = A_{D_j} Z(A_D) \]
which follows again from the universal property by localizing $A_{D_j}$ at the central element $D$. Hence, we can define a genuine non-commutative variety $\wis{spec}~\Ascr$ with central scheme $\wis{moduli}^{\theta}_{\alpha}~A$, finishing the proof of the main result and these talks.

\begingroup\raggedright\endgroup

\end{document}